\documentclass[review]{elsarticle} 

\usepackage{natbib}
\usepackage{lineno,hyperref}
\usepackage{amssymb,amsthm,latexsym, amsmath, amstext, amsfonts}
\usepackage[utf8]{inputenc}
\usepackage[english]{babel}
\usepackage{graphicx}
\usepackage{stackengine}
\usepackage{lineno}
\usepackage{xcolor}
\usepackage{bm}
\usepackage{changepage}

\usepackage{booktabs} 
\usepackage{colortbl}
\usepackage{xfrac}

\makeatletter
\def\ps@pprintTitle{%
   \let\@oddhead\@empty
   \let\@evenhead\@empty
   \def\@oddfoot{\reset@font\hfil\thepage\hfil}
   \let\@evenfoot\@oddfoot
}
\makeatother

\newcommand{\EqsModel}{\eqref{eq:orchard1_EqS1t}-\eqref{eq:orchard2_EqI2v}}
\newcommand{\bmh}[1]{ \bm{ \hat{ #1} } }
\newcommand{\deriv}[1]{ \dfrac{d #1}{dt} }

\newcommand{\comillas}[1]{``#1''}
\def\R{{\mathbf R}}

\modulolinenumbers[5]

\begin{document}

\begin{frontmatter}

\title{Investigating HLB control strategies using Genetic Algorithms: A two-orchard model approach with ACP Dispersal.}

\author[buap1]{A.~Anzo-Hern\'andez}
\ead{andres.anzo@hotmail.com}

\author[buap2]{U.J.~Giménez-Mujica}
\ead{uvencioj@gmail.com}

\author[buap1]{C.~Hernández-Gracidas}
\ead{cahernandezgr@conacyt.mx}

\author[buap2]{J.J.~Oliveros-Oliveros}
\ead{jacobo@gmail.com}

\address[buap1]{
{\textsc{CONHACYT - Benem\'erita Universidad Aut\'onoma de Puebla - Facultad de Ciencias F\'isico-Matem\'aticas},\\}
{Benem\'erita Universidad Aut\'onoma de Puebla},\\
\textsc{Avenida San Claudio y 18 Sur, Colonia San Manuel, 72570.\\
Puebla, Puebla, M\'exico.
\vskip 2ex}}

\address[buap2]{
{\textsc{Benem\'erita Universidad Aut\'onoma de Puebla - Facultad de Ciencias F\'isico-Matem\'aticas},\\}
\textsc{Avenida San Claudio y 18 Sur, Colonia San Manuel, 72570.\\
Puebla, Puebla, M\'exico.
\vskip 2ex}}

\setcounter{page}{1}

\begin{abstract}
This study focuses on the use of genetic algorithms to optimize control parameters in two potential strategies called mechanical and chemical control, for mitigating the spread of Huanglongbing (HLB) in citrus orchards. By developing a two-orchard model that incorporates the dispersal of the Asian Citrus Psyllid (ACP), the cost functions and objective function are explored to assess the effectiveness of the proposed control strategies. The mobility of ACP is also taken into account to capture the disease dynamics more realistically. Additionally, a mathematical expression for the global reproduction number ($R_{0}$) is derived, allowing for sensitivity analysis of the model parameters when ACP mobility is present. Furthermore, we mathematically express the cost function and efficiency of the strategy in terms of the final size and individual $R_{0}$ of each patch (i.e., when ACP mobility is absent). The results obtained through the genetic algorithms reveal optimal parameters for each control strategy, providing valuable insights for decision-making in implementing effective control measures against HLB in citrus orchards. This study highlights the importance of optimizing control parameters in disease management in agriculture and provides a solid foundation for future research in developing disease control strategies based on genetic algorithms.
\end{abstract}

\end{frontmatter}

\section{Introduction.}

Huanglongbing (HLB), also known as citrus greening, is one of the most devastating citrus diseases worldwide, with significant global impact on the citrus farming industry, causing substantial economic losses \cite{Li_2020, Singerman_2020, Hodges_2012}. It is caused by the bacterium \emph{Candidatus Liberibacter spp.} and transmitted by the Asian Citrus Psyllid (ACP), also named \emph{Diaphorina citri Kuwayama} \cite{Gottwald_2010}. HLB infection leads to severe consequences, including stunted growth, fruit deformation, yellow shoots on citrus trees and, ultimately, the death of infected trees \cite{Wang_2013}. 

Managing and controlling HLB presents a complex challenge, particularly in regions with multiple citrus orchard, where the disease can spread rapidly through ACP dispersal. Effective management strategies for HLB disease encompass a range of measures. These include thorough inspection and removal of infected plants, application of targeted insecticide sprays to reduce the population of ACP vectors \cite{Boina_2015}, biocontrol agents\cite{Milosavljevic_2017}, utilization of healthy nursery trees for replanting, implementing differentiated psyllid control strategies based on tree location and shoot flushing patterns, coordinating area-wide psyllid control efforts, eliminating potential sources of infection in noncommercial properties near commercial orchards, among others. \cite{Alquezar_2022, Bassanezi_2020}. 

Identifying effective control strategies among numerous potential interventions can be challenging and costly to test in the field. A more efficient approach would be to prioritize interventions that show the highest potential for success in reducing the disease rate or severity, while considering the costs associated with implementing the strategy. The question then arises: how can we identify these strategies? Mathematical models offer a valuable tool for evaluating the cost-effectiveness of individual and combined intervention strategies in a timely manner. By providing insights into the effectiveness of different interventions, these models can help streamline experimental design and guide data collection efforts towards critical factors and interventions that are most likely to yield positive outcomes.

An example of the use of mathematical modeling to study control strategies for HLB disease was demonstrated in the work by C. Chiyaka \emph{et al.} \cite{Chiyaka_2012}. Their compartmental model highlighted that the effectiveness of psyllid spraying is influenced by factors such as the timing of initial spraying, frequency of application, and the efficacy of the insecticides used. Additionally, the impact of removing symptomatic flush was shown to be dependent on the frequency of removal and the timing of initiation relative to the start of the epidemic. Furthermore, J.A. Lee \emph{et al.} \cite{Lee_2015} introduced spatially explicit models that combined compartmental and agent-based approaches to simulate the transmission of HLB disease within a specific citrus orchard. Their research demonstrated that implementing intervention strategies aimed at reducing the psyllid population by 75\% during the flushing periods can effectively delay the infection of the entire grove. This delay in infection not only mitigates the spread of the disease, but also leads to a reduction in the overall amount of insecticide required throughout the year.

Control theory has also been employed to explore alternative control methods for managing HLB disease in dynamic systems. For example, W. Ling \emph{et al.} contributed to the development of control strategies for HLB disease by utilizing mathematical models based on differential equations with discontinuous right-hand sides \cite{Ling_2021}. These models incorporated the discontinuous removal of infected trees as a control strategy for managing the spread of HLB disease. The application of optimal control techniques to manage HLB transmission within orchards has gained significant attention in the literature. Researchers such as F. Zhang \emph{et al.} \cite{FZhang_2020, FZhang_2021}, R. G. Vilamiu \cite{Vilamiu_2012}, and R. Taylor \cite{Taylor_2016} have utilized these methods to investigate the effectiveness of various management strategies. By considering the cost of implementing different control measures as a profit measure, these studies aimed to identify the optimal strategy that maximizes outcomes while minimizing costs. Such an approach provides valuable insights for designing efficient and cost-effective HLB management strategies in orchards.

Despite the significant progress made in mathematical modeling and control theory for studying HLB propagation, the current research primarily focuses on the spread of the bacteria within individual orchards, without considering the emigration of the ACP, which plays a significant role in the propagation of HLB. An empirical study conducted by H. Lewis-Rosenblum \emph{et. al.} in \cite{Lewis_Rosenblum_2015} showed the long-range dispersal capabilities of the ACP, and revealed that it is capable of traversing potential geographical barriers, covering distances of at least 2 km within a 12-day period. Furthermore, the presence of abandoned orchards \cite{Tiwari_2010} and other ecosystems such as lakes \cite{Martini_2013} significantly affects the ability of ACP to undertake long-distance travel. These ecosystems serve as temporary habitats for ACP, and when favorable conditions arise, they migrate to more suitable locations. Notably, relative humidity, temperature, and light periods influence ACP emigration \cite{Stelinski_2019, Tomaseto_2017}, as well as seasonality \cite{Zorzenon_2020, Hall_2011}. Even insecticide applications can impact ACP emigration \cite{Johnston_2019}.

In this paper, we propose an epidemiological mathematical model for two orchards, considering the mobility of ACP between them. The model is based on the concept of metapopulation networks, which have gained interest in epidemiology for studying human mobility between cities or urban areas \cite{Uvencio_2022}. Drawing inspiration from previous studies on Dengue disease \cite{Anzo_2019}, our model incorporates two populations: citrus trees and ACP. Unlike the previous model, where mosquitoes are stationary and humans are mobile, our model assumes that ACP can move while citrus trees remain static. We also introduce two forms of control parameters: mechanical and chemical control. Mechanical control involves using physical barriers, such as yellow-colored traps or cultural practices, to enhance orchard vigilance. In contrast, chemical control involves fumigation campaigns with insecticides applied throughout the entire orchard. To optimize the cost and effectiveness of these control strategies, we employ genetic algorithms. These algorithms, inspired by the process of natural selection, involve the generation of diverse sets of potential solutions and iteratively refining them to find the most optimal parameters. This iterative process mimics the way in which nature evolves over time, enabling us to systematically explore a wide range of control parameter combinations. By doing so, we can identify those combinations that lead to the desired outcomes with minimal costs. This approach not only enhances our ability to achieve effective disease control but also provides a systematic framework for decision-making in the complex domain of agricultural disease management.

\section{Preliminaries.}

The aim in this section is to describe the proposed mathematical model for HLB transmission between two citrus orchards with parametric control, which is based on the compartmental SAIR$-$SI model. The local basic reproduction number $R_{0}$ is calculated for this model through the application of the Next Generation Matrix (NGM) approach.

\subsection{Model description.}

The mobility of the ACP between orchards plays a crucial role in the spread of HLB infections. In order to develop a mathematical model for HLB transmission between two orchards, we identified two potential sources of bacterial transmission within a single orchard: infected ACP present in their own orchard, or those that have traveled from neighboring orchards. Similarly, susceptible ACP can get the infection in two ways: by feeding on infected trees in their own orchards or by traveling to a neighboring orchard and becoming infected there. 

To incorporate ACP dispersal into a compartmental model for HBL transmission, we first assume that the dynamics of bacteria transmission among citrus trees in each orchard is driven by the SAIR-SI model. This model divides, at any time instant $t$, the epidemic state of the citrus trees in orchard with label $i$ ($=1,2$) into four categories: Susceptible trees $S_{i\tau}(t)$, infectious and asymptomatic trees $A_{i\tau}(t)$ (trees that have acquired the bacteria but do not display symptoms), infectious and symptomatic trees $I_{i\tau}(t)$ (citrus trees that display symptoms), and roguing trees $R_{i\tau}(t)$ (trees that have been cut down and removed from the orchard by the farmer). Also, the epidemic state for ACP vectors is divided into susceptible and infectious psyllids $S_{iv}(t)$ and $I_{iv}(t)$, respectively. From here onwards, the first subscript $i=1,2$ denotes the patch label in both the variable states and the parameters. 

We introduce the parameter $\phi_{12} \in [0,1]$ to represent the average fraction of ACP from orchard 1 that is present in orchard 2 at any given time, thereby accounting for emigration rates from orchard 1 to orchard 2. It is important to note that in our model, we assume that ACP is capable of traveling over greater distances with a specific direction of dispersal. For this reason, we set $\phi_{21} = 0$, indicating that ACP vectors can only travel from orchard 1 to orchard 2 and not vice versa. This assumption reflects the directional mechanism of ACP dispersal in our model.

In this context, the rates of immigration for susceptible and infectious ACP coming from orchard 1 and traveling to orchard 2 are quantified by $\phi_{12}S_{v1}$ and $\phi_{12}I_{v1}$, respectively. Additionally, we use $\phi_{11} = 1-\phi_{12}$ to represent the fraction of ACP that remain in their own orchard during HLB transmission. With the inclusion of ACP dispersal, the two-orchard SAIR-SI model with parametric control can be represented by a system of twelve differential equations. Specifically, six of these equations capture the HLB disease transmission dynamics in orchard 1 according to:

\begin{eqnarray} 
\deriv{S_{1\tau}} &=&  -\bmh{\beta}_{1\tau} \dfrac{S_{1\tau}(\phi_{11}I_{1v})}{N_{1\tau}} \, , \label{eq:orchard1_EqS1t} \\ 
\nonumber \\
\deriv{A_{1\tau}} &=&  \bmh{\beta}_{1\tau} \dfrac{S_{1\tau}(\phi_{11}I_{1v})}{N_{1\tau}} - (\sigma_{1} + \mu_{1\tau})A_{1\tau} \, ,  \label{eq:orchard1_EqA1t} \\ 
\nonumber \\ 
\deriv{I_{1\tau}} &=&  \sigma_{1} A_{1\tau} - \bmh{r}_{1\tau} I_{1\tau}  \, ,  \label{eq:orchard1_EqI1t} \\ 
\nonumber \\
\deriv{R_{1\tau}} &=&  \mu_{1\tau}A_{1\tau} + \bmh{r}_{1\tau}I_{1\tau} \, ,     \label{eq:orchard1_EqR1t} \\ 
\nonumber \\ 
\deriv{S_{1v} } &=& \Lambda_{1} - \bmh{\beta}_{1v} \dfrac{\phi_{11}S_{1v} (A_{1\tau} + I_{1\tau})}{N_{1\tau}} -  \bmh{\beta}_{2v} \dfrac{\phi_{12}S_{1v} (A_{2\tau} + I_{2\tau})}{N_{2\tau}} - \bmh{\mu}_{1v}S_{1v}, \nonumber \\ 
\label{eq:orchard1_EqS1v}  \\ 
\deriv{I_{1v} }  &=&  \bmh{\beta}_{1v} \dfrac{\phi_{11}S_{1v} (A_{1\tau} + I_{1\tau})}{N_{1\tau}} + \bmh{\beta}_{2v} \dfrac{\phi_{12}S_{1v} (A_{2\tau} + I_{2\tau})}{N_{2\tau}} - \bmh{\mu}_{1v}I_{1v} \, .  \label{eq:orchard1_EqI1v}
\end{eqnarray}

\begin{figure}[ht]
\begin{center}
\includegraphics[width=12.0cm,height=9cm]{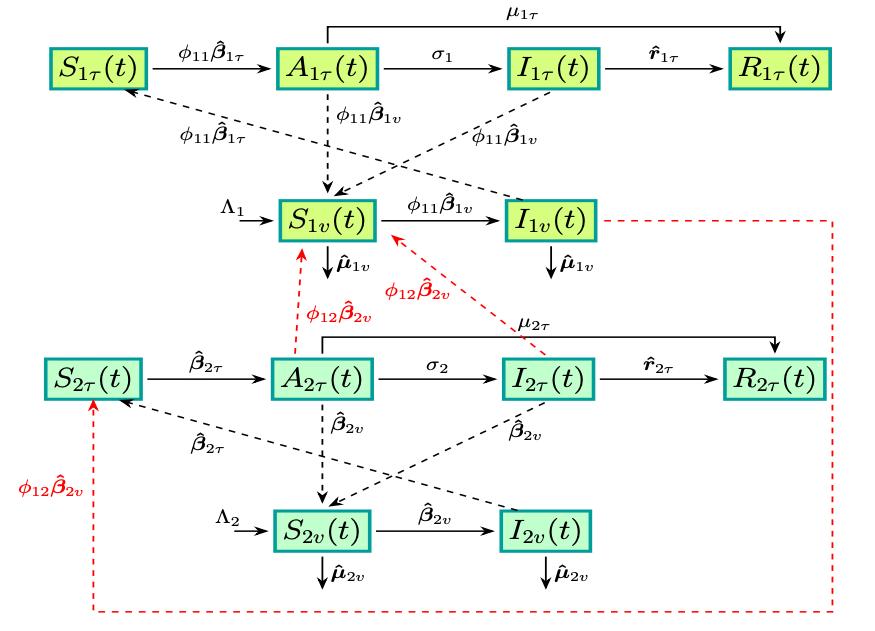}  
\caption{Block diagram representing the two-orchards SAIR-SI model with parametric control and ACP dispersal (Eqs. \EqsModel).} 
\label{fig:HLB_SAIR_SI_Diagram}
\end{center}
\end{figure}

For orchard 2, the HLB disease transmission dynamics is described by the following set of six equations: 
\begin{eqnarray} 
\deriv{S_{2\tau}}  &=&  - \bmh{\beta}_{2\tau}\dfrac{S_{2\tau} I_{2v}}{N_{2\tau}} - \bmh{\beta}_{2\tau}\dfrac{S_{2\tau}(\phi_{12}I_{1v})}{N_{2\tau}} \, ,  \label{eq:orchard2_EqS2t} \\
\nonumber \\
\deriv{A_{2\tau}}  &=&  \bmh{\beta}_{2\tau}\dfrac{S_{2\tau} I_{2v} }{N_{2\tau}}  + \bmh{\beta}_{2\tau}\dfrac{S_{2\tau}(\phi_{12}I_{1v}) }{N_{2\tau}} - (\sigma_{2} + \mu_{2\tau})A_{2\tau} \, ,  \label{eq:orchard2_EqA2t}   \\ 
 \nonumber \\
\deriv{I_{2\tau}} &=&  \sigma_{2} A_{2\tau} - \bmh{r}_{2\tau} I_{2\tau}  \, , \label{eq:orchard2_EqI2t} \\
\nonumber \\
\deriv{R_{2\tau}} &=&  \mu_{2\tau}A_{2\tau} + \bmh{r}_{2\tau}I_{2\tau} \, ; \label{eq:orchard2_EqR2t} \\
\nonumber \\
\deriv{S_{2v} }   &=& \Lambda_{2}  - \bmh{\beta}_{2v} \dfrac{S_{2v} (A_{2\tau} + I_{2\tau})}{N_{2\tau}} - \bmh{\mu}_{2v}S_{2v} \, ,  \label{eq:orchard2_EqS2v} \\
 \nonumber\\
\deriv{I_{2v}}   &=&  \bmh{\beta}_{2v} \dfrac{S_{2v}(A_{2\tau} + I_{2\tau})}{N_{2\tau}} - \bmh{\mu}_{2v}I_{2v} \, ;  \label{eq:orchard2_EqI2v}
\end{eqnarray}

\noindent where $\mu_{i\tau}$ denotes the mortality rate of citrus trees, $\sigma_{i}$ represents the progression rate from asymptomatic to symptomatic infection in trees and $\Lambda_{i}$ is the constant recruitment rate of ACP vectors, with $i=1,2$ being the orchard label. The recruitment rate parameter captures the processes of reproduction and development of ACP vectors within the orchard, including factors such as the rate at which ACP vectors reproduce, the number of eggs laid per female (fecundity), and the time it takes for the eggs to develop into adults. Figure \eqref{fig:HLB_SAIR_SI_Diagram} illustrates a block diagram representing the two-orchard SAIR-SI model with parametric control and ACP dispersal.

In each orchard, we have incorporated four control parameters, denoted as $\bmh{\beta}_{i\tau}$, $\bmh{r}_{i\tau}$, $\bmh{\beta}_{iv}$, and $\bmh{\mu}_{iv}$. In Eqs. \EqsModel, we have used the hat symbol and bold typography to simplify the notation of the control parameters. These parameters are associated with two distinct types of control strategies:

\begin{table}[ht]
\LARGE
\centering\resizebox{13.0 cm}{!}{
\begin{tabular}{clccr}
\toprule 
\textbf{Parameter} & \hspace{3cm} \textbf{Description} & \textbf{Baseline value} &\textbf{Range} & \textbf{Ref.} \\
\midrule 
$N_{i\tau}$ & Total number of citrus trees in the orchard $i$ & 2000 & [1000,2000]  & \cite{FZhang_2020, Luo_2017} \\
\rowcolor{black!20}
$\mu_{i\tau}$ & Natural mortality rate of citrus trees ($month^{-1}$) & 0.0033 & [0.00275,0.00417]  & \cite{FZhang_2020, Wang_2014} \\ 
$\mu_{iv}$ & Natural mortality rate of ACP ($month^{-1}$) & 0.266 & [0.117,0.95]  & \cite{Liu_2000, FZhang_2020} \\ 
\rowcolor{black!20}
$\pi_{\tau}$ & Probability of HLB transmission from infected ACP to susceptible trees & 0.025 & [0,1) & \cite{Taylor_2016,Vilamiu_2012} \\
$\pi_{v}$ & Probability of HLB transmission from infected tree to susceptible ACP & 0.13 & [0,1) & \cite{Taylor_2016,Vilamiu_2012} \\
\rowcolor{black!20}
$b$       & Probing rate of tress & 1.15  & [1,1.15] & \cite{Taylor_2016,Vilamiu_2012} \\
$r_{i\tau}$   & Effectiveness of orchard vigilance (Number of roguing trees per month)   &  0.47 & [0,1)  & \cite{Vilamiu_2012} \\
\rowcolor{black!20}
$\sigma_{i}$  & Rate at which asymptomatic infected citrus tree become symptomatic ($month^{-1}$). & 0.2 & 0.155-0.99 & \cite{Chiyaka_2012,FZhang_2020}\\
$\omega_{i}$    & Max abundance of ACP vectors per citrus tree & 100 & 50-100   &  \cite{Wang_2014} \\
\rowcolor{black!20}
$\Lambda_{i}$ & Recruitment rate of ACP ($month^{-1}$) & 25,000 & [25,000,50,000] & \cite{FZhang_2020, Sule_2012} \\
$\phi_{12}$    & Average fraction of ACP vectors from orchard 1 that are present in orchard 2 & 0.35 & [0,1)   & Assumed \\ 
\rowcolor{black!20}
$m_{i}, n_{i}$    & Control parameter for a mechanical strategy & -- & [0,1) & Assumed \\
$p_{i}, q_{i}$    & Control parameter for a chemical strategy   & -- & [0,1) & Assumed \\
\bottomrule
\end{tabular} }
\caption{Description and parameter values for the two-orchard SAIR$-$SI model (Eqs. \EqsModel).}
\label{table_parameters}
\end{table}

\textbf{Mechanical strategy}: This strategy aims to reduce the interaction between citrus trees and ACP vectors. In order to implement this strategy, we focus on modifying two key parameters:

\begin{itemize}
    \item[$\ast$] $\bmh{\beta}_{i\tau} = (1-m_{i})b \pi_{\tau}$ which represents the rate at which susceptible citrus plants become infected through contact with infectious ACP vectors. Here, $\pi_{\tau}$ is a fixed parameter indicating the probability of HLB transmission from ACP vector to trees, $b_{i}$ represents the probing rate of host trees, that is, it is a measure of the rate at which the insect feeds on citrus trees and transmits the pathogen (this parameter reflects the vector's behavior and its feeding activity); and  $m_{i} \in [0,1)$ is the control parameter whose aim is to reduce the value of $b_{i}$. Some of the control strategies that can be employed to control $\bmh{\beta}_{i\tau}$ are:  1) Protective barriers: Physical barriers can be used to physically prevent psyllids from reaching the citrus trees. These barriers can include fine-mesh nets or screens installed around the trees or the entire orchard to create a physical barrier between the vectors and the trees. 2) Trapping and monitoring: Placing sticky traps or yellow-colored traps in the orchard can help capture and monitor ACP populations. 3) Cultural practices: Implementing cultural practices such as pruning, weed control, and orchard sanitation can help reduce the habitat and food sources for psyllids, thereby reducing their probing and feeding activities.
    
    \item [$\ast$] $\bmh{r}_{i\tau} = 1-n_{i}$ quantifies the level of orchard vigilance and the effectiveness of identifying trees with clear symptomatic signals, such as the characteristic mottled pattern of yellowing in leaves, misshapen fruits, and bitterness. The control parameter $n_{i} \in [0,1)$ aims to adjust the orchard vigilance  parameter $\bmh{r}_{i\tau}$ in a range of values between zero and one, where values close to zero indicate few vigilance and close to one a daily and strict vigilance. In a study by Vilamiou \emph{et al.} (\cite{Vilamiu_2012}), they implemented a control strategy involving the removal of symptomatic trees every 3 months, with an efficacy of 47\% ($r_{i\tau} = 0.47$). Several control strategies can be employed to manage $\bmh{r}_{i\tau}$, including: 1) Training and education: Provide comprehensive training and education to farmers, orchard workers, and inspectors regarding the symptoms of HLB disease in citrus trees. 2) Regular scouting and monitoring: Establish a systematic program for regularly inspecting trees in the orchard, focusing on identifying symptoms associated with HLB disease. 3) Rapid diagnostic tools: Utilize advanced diagnostic tools, such as PCR tests, immunological assays, or DNA-based techniques, to enable early detection of HLB disease.  
\end{itemize}

\textbf{Chemical strategy}: This strategy aims to decrease the population of ACP vectors and increase their mortality rate. To effectively implement this strategy, we focus on modifying two crucial parameters:

\begin{itemize}
    \item[$\ast$] $\bmh{\beta}_{iv} = (1-p_{i})\omega_{i} \pi_{v}$ which represents the rate at which susceptible ACP vectors become infected, where $\pi_{v}$ is a fixed parameter indicating the probability of HLB transmission from infectious critrus tree to susceptible ACP vectors; $\omega_{i}$ is the maximum abundance of ACP vectors per citrus tree and $p_{i} \in [0,1)$ is the control parameter whose aim is to reduce the value of $\omega_{i}$. Some of the control strategies that can be employed to control $\bmh{\beta}_{iv}$ are: 1) Trapping and Removal: Deploy traps designed to attract and capture ACP vectors. These traps can use pheromones or visual cues to lure the insects. Once captured, the vectors can be removed from the orchard, reducing their population. 2) Biological Control: Implement biological control methods by introducing natural enemies of ACP vectors, such as parasitic wasps or predatory insects. These natural enemies can help suppress the population of ACP vectors without relying heavily on chemical insecticides. 3) Cultural Practices: Adopt cultural practices that discourage ACP vectors, such as removing weeds or alternate host plants that serve as reservoirs for the vectors. Proper sanitation and pruning techniques can also help reduce their population.
    
    \item [$\ast$] $\bmh{\mu}_{iv} = (1+q_{i})\mu_{iv}$, where $\mu_{iv}$ is the baseline value of the ACP mortality rate. The objective of this control strategy is to increase the mortality of ACP vectors in the orchard $i$. Typically, this involves implementing fumigation campaigns where insecticides are sprayed throughout the entire orchard. The effectiveness of this strategy can be achieved by either increasing the value of $\mu_{iv}$, which directly influences the mortality rate, or by reducing the recruitment rate parameter $\Lambda_{i}$. However, for the sake of model simplicity, we do not consider $\Lambda_{i}$ as a control parameter in our model.  Some of the control strategies that can be employed to control $\bmh{\mu}_{iv}$ are: 1) apply insecticides specifically targeted at ACP vectors to increase their mortality rate, or 2) use systemic insecticides that are absorbed by the citrus trees and distributed throughout their tissues. When ACP vectors feed on the trees, they ingest the insecticide, which can help control their population.
\end{itemize}

To simplify the complexity arising from the number of control parameters, we propose categorizing the SAIR-SI model into two scenarios. The first scenario is the two-orchard SAIR-SI model with mechanical control, where the chemical control parameters are set to $p_{i} = q_{i} = 0$ for $i=1, 2$. In this scenario, we focus solely on mechanical control methods. The second scenario is the two-orchard SAIR-SI model with chemical control, where chemical control is applied in the orchards and the mechanical control parameters are set $m_{i} = n_{i} = 0$. In this scenario, we consider the use of chemical control methods.

We present a summary of the parameter descriptions for two-orchard SAIR-SI model presented in Eqs. \EqsModel in Table \eqref{table_parameters}, which includes numerical values reported in references for easy comparison.

\begin{figure}[t!]
\begin{center}
\includegraphics[width=12.5cm,height=8cm]{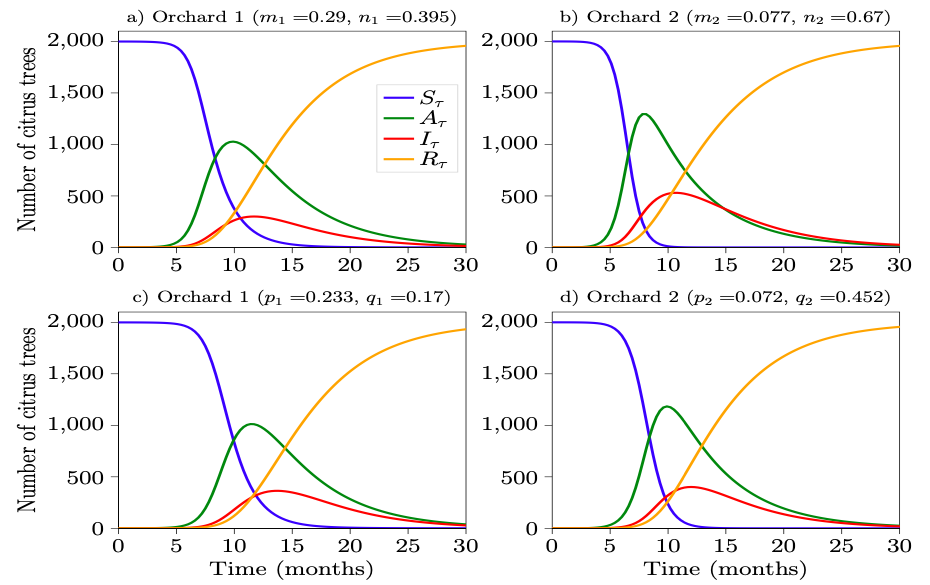}
\caption{Numerical solution to the two orchard model for HLB transmission with $\phi_{12} = 0.35$ and a) mechanical control with $m_{1} = 0.2$ and  $n_{1} = 0.6$ for orchard 1, and b) $m_{2} = 0.47$ and $n_{2} = 0.6$ for orchard 2; c) chemical control with  $p_{1} = 0.45$ and $q_{1} = 0.24$ for orchard 1 and d) $p_{2} = 0.57$ and $q_{2} = 0.8$ for orchard 2.} 
\label{fig:SAIR_SI_Sol}
\end{center}
\end{figure}

To illustrate the dynamics of HLB propagation between two orchards, we present the numerical solution to Eqs.\EqsModel with mechanical and chemical control strategies in Figure \ref{fig:SAIR_SI_Sol}. For the mechanical control strategy, we set the values of $m_{1} = 0.29$ and $n_{1} = 0.395$ for orchard 1 and; $m_{2} = 0.07$ and $n_{2} = 0.67$ for orchard 2. Similarly, for the chemical control strategy, we set $p_{1} = 0.23$ and $q_{1} = 0.17$ for orchard 1 and; $p_{2} = 0.072$ and $q_{2} = 0.45$ for orchard 2, denoting the respective control levels.

Considering the baseline values for the rest of the parameters described in Table \eqref{table_parameters}, the initial conditions for orchard 1 are set as:
\begin{align}\label{eq:init_cond_orchar_1}
S_{1\tau}(t=0) &= N_{1\tau}, \quad A_{1\tau}(t=0) = 1, \quad I_{1\tau}(t=0) = 0, \nonumber \\
R_{1\tau}(t=0) &= 0, \quad S_{1v}(t=0) = \Lambda_{1}/\mu_{1v}, \quad \text{and} \quad I_{1v}(t=0) = 0 \,.
\end{align}

\noindent For orchard 2, we assume the absence of infectious trees or infectious ACP, and accordingly, we set its initial conditions as follows:
\begin{align}\label{eq:init_cond_orchar_2}
S_{2\tau}(t=0) &= N_{2\tau}, \quad A_{2\tau}(t=0) = 0, \quad I_{2\tau}(t=0) = 0, \nonumber \\
R_{2\tau}(t=0) &= 0, \quad S_{2v}(t=0) = \Lambda_{2}/\mu_{2v}, \quad \text{and} \quad I_{2v}(t=0) = 0 \,;
\end{align}

We set the dispersal parameter at $\phi_{12} = 0.35$,representing the assumption that approximately 35\% of the ACP population in orchard 1 migrates to the disease-free orchard 2. This choice is made to provide a baseline scenario for studying the impact of ACP mobility in our model.

Based on the findings from the numerical example, as depicted in Figure \eqref{fig:SAIR_SI_Sol}, it is evident that even a small proportion of ACP dispersal can lead to the transmission of HLB bacteria to a neighboring orchard, regardless of the implemented control measures. Moreover, it is important to mention that the highest point of the infectious asymptomatic and symptomatic curves in orchard 2 may occur ten or twelve months after the initial onset of the HLB epidemic. Moreover, we observe that the implementation of the mechanical control strategy results in a maximum of $1,027$ infectious asymptomatic trees and $300$ symptomatic trees in orchard 1. Similarly, orchard 2 experiences a maximum of $1,294$ asymptomatic trees and $530$ symptomatic trees. On the other hand, the employment of the chemical control strategy leads to a maximum of $1,012$ infectious asymptomatic trees and $364$ symptomatic trees in orchard 1. In orchard 2, these numbers reach 1,183.79 asymptomatic trees and 401.528 symptomatic trees.

\subsection{Basic reproduction number and its parametric sensitivity analysis.}

The basic reproduction number $R_{0}$ is used in epidemiology to quantify the transmission potential of an infectious disease in a population. In the context of HLB transmission, it quantifies the secondary infections caused by a single infected citrus tree in a given orchard, and how it may be influenced by ACP dispersal to neighboring orchards. In this section, we utilize the Next Generation Matrix (NGM) methodology to mathematically express the global reproduction number for two orchards, referred to as $R_{\text{g0}}$, in terms of the individual orchard reproduction numbers $R_{10}$ and $R_{20}$ referred to the basic reproduction number of the individual and isolated orchards when ACP dispersal is absent (i.e., $p_{12} = 0$).

The NGM method uses the infection as a demographic process where newly infected individuals are added or removed \cite{Diekmann_2010}. The initial step in the NGM method involves dividing the compartmental model into two sub-systems: the infected sub-system, also referred to as the infected compartments, and the disease-free sub-system. For the two-orchard SAIR$-$SI model with parametric control, the infected sub-system comprises the state variables $x = (A_{1\tau}, I_{1\tau}, I_{1v}, A_{2\tau}, I_{2\tau}, I_{2v}) \in \R^{6}$, while the state variables $y = (S_{1\tau}, R_{1\tau}, S_{1v}, S_{2\tau}, R_{2\tau}, S_{2v}) \in \R^{6}$ make up the disease-free sub-system. The right-hand side of the infected sub-system is then separated as follows:

\[
\begin{array}{lll}
\deriv{A_{1\tau}} &=&  \underbrace{ \bmh{\beta}_{1\tau}\dfrac{S_{1\tau}( \phi_{11}I_{1v} ) } {N_{1\tau}}}_{\mathbf{F}_{11}}  \quad \underbrace{  - (\sigma_{1} + \mu_{1\tau})A_{1\tau}  }_{  \mathbf{V}_{11}} \, ,  \\
\\
\deriv{I_{1\tau}}    &=&  \underbrace{\sigma_{1} A_{1\tau} - \bmh{r}_{1\tau} I_{1\tau}}_{\mathbf{V}_{12}}  \, ,  \\
\\
\deriv{I_{1v}}       &=&   \underbrace{\bmh{\beta}_{1v} \dfrac{\phi_{11}S_{1v}(A_{1\tau} + I_{1\tau})}{N_{1\tau}} + \bmh{\beta}_{2v} \dfrac{\phi_{12}S_{1v}(A_{2\tau} + I_{2\tau})}{N_{2\tau}} }_{\mathbf{F}_{13}}  \underbrace{- \bmh{\mu}_{1v}I_{1v}}_{\mathbf{V}_{13}} \, ;    \\
\end{array}
\]

\noindent and
\[
\begin{array}{lll}
\deriv{A_{2\tau}}     &=&  \underbrace{\bmh{\beta}_{2\tau}\dfrac{S_{2\tau}I_{2v}} {N_{2\tau}}  + \bmh{\beta}_{2\tau}\dfrac{S_{2\tau}(\phi_{12}I_{1v})}{N_{2\tau}}}_{\mathbf{F}_{21}}  \quad \underbrace{  - (\sigma_{2} + \mu_{2\tau})A_{2\tau}  }_{  \mathbf{V}_{21}} \, ,  \\
\\
\deriv{I_{2\tau}}    &=&  \underbrace{\sigma_{2} A_{2\tau} - \bmh{r}_{2\tau} I_{2\tau}}_{\mathbf{V}_{22}}  \, ,  \\
\\
\deriv{I_{2v}}       &=&   \underbrace{ \bmh{\beta}_{2v} \dfrac{S_{2v}(A_{2\tau} + I_{2\tau})}{N_{2\tau}}}_{\mathbf{F}_{23}}  \underbrace{- \bmh{\mu}_{2v}I_{2v}}_{\mathbf{V}_{23}} \, ;    \\
\end{array}
\]

\noindent where $\mathbf{F}_{ij}$ are the transition terms and $\mathbf{V}_{ij}$ the birth, death, disease progression or recovery terms, for $i,j=1,2$. 

Next, the NGM methodology involves constructing the following matrices:
\[
\mathcal{F} = \left. \left[ \scalebox{1.5}{$\frac{\partial \mathbf{F}_{ij}}{\partial x_{k}}$} \right|_{\tiny{P}} \right]_{6 \times 6} \quad \text{and} \quad 
\mathcal{V} = \left. \left[ \scalebox{1.5}{$\frac{\partial \mathbf{V}_{ij}}{\partial x_{k}}$} \right|_{\tiny{P}} \right]_{6 \times 6} \, , 
\]

\noindent where  $P = (N_{1\tau}, 0, 0, 0, \Lambda_{1}/\mu_{1v}, 0,  N_{2\tau}, 0, 0, 0, \Lambda_{2}/\mu_{2v}, 0) \in \R^{12}$ represents the disease-free equilibrium state (\emph{DFE}). Then:
\[
\mathcal{F} =
\begin{pmatrix}
0 & 0 & \bmh{\beta}_{1\tau}\phi_{11} & 0 & 0 & 0 \\
0 & 0 & 0               & 0 & 0 & 0\\
\bmh{\beta}_{1v}\phi_{11}\dfrac{N^{0}_{1v}}{N_{1\tau}} &  \bmh{\beta}_{1v}\dfrac{N^{0}_{1v}}{N_{1\tau}} & 0 & \bmh{\beta}_{2v}\phi_{12}\dfrac{N^{0}_{1v}}{N_{2\tau}} & \bmh{\beta}_{2v}\phi_{12}\dfrac{N^{0}_{1v}}{N_{2\tau}} & 0 \\
0           & 0 & \bmh{\beta}_{2\tau}\phi_{12} & 0 & 0 & \bmh{\beta}_{2\tau} \\ 
0           & 0 & 0                & 0 & 0 &  0         \\
0           & 0 & 0                & \bmh{\beta}_{2v}\dfrac{N^{0}_{2v}}{N_{2\tau}} & \bmh{\beta}_{2v}\dfrac{N^{0}_{2v}}{N_{2\tau}} & 0
\end{pmatrix}; 
\]

\noindent and
\[
\mathcal{V} = 
\begin{pmatrix}
-\sigma_{1} -\mu_{1\tau} &  0               & 0          & 0 & 0 & 0 \\
\sigma_{1}               & -\bmh{r}_{1\tau} & 0          & 0 & 0 & 0 \\
0                        & 0                & -\bmh{\mu}_{1v}   & 0 & 0 & 0  \\
0                        & 0                & 0          & -\sigma_{2} -\mu_{2\tau} & 0 & 0 \\
0                        & 0                & 0          & \sigma_{2}               & -\bmh{r}_{2\tau}   & 0 \\
0                        & 0                & 0   & 0 & 0 & -\bmh{\mu}_{2v}
\end{pmatrix} \, ;
\]

\noindent where $N^{0}_{iv} = \Lambda_{i}/\mu_{iv}$. The matrix $K = -\mathcal{F}\mathcal{V}^{-1}$ is called the NGM \cite{Diekmann_2010} which, for the two-orchard model with parametric control, becomes: 

\begin{adjustwidth}{-0.5cm}{-0.5cm} 
\[
K = \begin{pmatrix}
    0   & 0 & \dfrac{\bmh{\beta}_{1\tau}\phi_{11}}{\bmh{\mu}_{1v}} & 0 & 0 & 0 \\
    0   & 0 &  0                                & 0 & 0 & 0 \\
    \bmh{\beta}_{1v}\dfrac{N^{0}_{1v}}{N_{1\tau}}\dfrac{(\sigma_{1} + \phi_{11}\bmh{r}_{1\tau})}{\bmh{r}_{1\tau}(\sigma_{1} + \mu_{1\tau})} & \bmh{\beta}_{1v}\dfrac{N^{0}_{1v}}{\bmh{r}_{1\tau}N_{1\tau}} & 0 & \bmh{\beta}_{2v}\phi_{12}\dfrac{N^{0}_{1v}}{N_{2\tau}}\dfrac{(\sigma_{2} + \bmh{r}_{2\tau})}{\bmh{r}_{2\tau}(\sigma_{2} + \mu_{2\tau})} & \bmh{\beta}_{2v}\phi_{12}\dfrac{N^{0}_{1v}}{\bmh{r}_{2\tau}N_{2\tau}} & 0 \\
    0 & 0 & \dfrac{\bmh{\beta}_{2\tau}\phi_{12}}{\bmh{\mu}_{1v}} & 0 & 0 & \dfrac{\bmh{\beta}_{2\tau}}{\bmh{\mu}_{2v}}\\
    0   & 0 &  0                                & 0 & 0 & 0 \\
    0   & 0 &  0 & \bmh{\beta}_{2v}\dfrac{N^{0}_{2v}}{N_{2\tau}}\dfrac{(\sigma_{2} + \bmh{r}_{2\tau})}{\bmh{r}_{2\tau}(\sigma_{2} + \mu_{2\tau})} & \bmh{\beta}_{2v}\dfrac{N^{0}_{2v}}{\bmh{r}_{2\tau}N_{2\tau}} & 0 
    \end{pmatrix}.
\]
\end{adjustwidth}

Therefore, the two-orchard basic reproduction number $R_{\text{0g}}$ is defined as the spectral radius $\rho(K)$ of the NGM. Then, the characteristic polynomial $p(\lambda) = det( K - \lambda \mathbb{I}_{6}) = 0$, where $\mathbb{I}_{6}$ is the $6 \times 6$ identity matrix given by:
\begin{equation}\label{eq:characteristic_polynomial}
    \lambda^{6} - \underbrace{(R^{2}_{10} + (1+\phi_{12}\delta)\cdot R^{2}_{20} )}_{b} \lambda^{4} - \underbrace{R^{2}_{10}\cdot R^{2}_{20}}_{c}\lambda^{2} = 0
\end{equation}

\noindent where $\delta = \delta_{1}/\delta_{2}$, with $\delta_{1}= N^{0}_{1v}/\bmh{\mu}_{1v}$ and $\delta_{2}= N^{0}_{2v}/\bmh{\mu}_{2v}$ being the population growth rate of ACP in orchard one and two, respectively. The terms $R_{10}$ and $R_{20}$ represents the corresponding local basic reproduction number for orchard one and two when ACP dispersal is not considered, that is, when $p_{12} = 0$. Following the NGM methodology described in this section, it is straightforward to show that: 
\begin{equation}\label{eq:Rio}
\begin{array}{lll}
R_{i0} & = & \sqrt{\dfrac{\bmh{\beta}_{i\tau}}{\sigma_{i} + \mu_{i\tau}} \cdot \dfrac{\bmh{\beta}_{iv}}{\bmh{\mu}_{iv}} \cdot \dfrac{N^{0}_{iv}}{N_{i\tau}} \cdot \Big( \dfrac{\sigma_{i}}{\bmh{r}_{i\tau}} + 1\Big)}  \, ,\\
\\
& = & \sqrt{R_{i}^{\tau v} \cdot R_{i}^{v \tau} \cdot \dfrac{\Lambda_{i}}{\mu_{iv}N_{i\tau}} \cdot \Big( \dfrac{\sigma_{i}}{\bmh{r}_{i\tau} + 1} \Big) }\, ;
\end{array}
\end{equation}

\noindent for $i = 1,2$; where $R_{i}^{\tau v} = \bmh{\beta}_{i\tau}/(\sigma_{i} + \mu_{i\tau})$ denotes the ratio of new infections generated by a single asymptomatic infectious tree during the period of time it belongs in the orchard, \emph{i.e.} $1/(\sigma_{i} + \mu_{i\tau})$. Analogously, $R^{v\tau}_{i} = \bmh{\beta}_{i2}/\bmh{\mu}_{v}$ is the ratio of new infections caused by an ACP vector during its lifespan, \emph{i.e.} $1/\bmh{\mu}_{v}$. 

Then, by solving the characteristic polynomial Eq. \eqref{eq:characteristic_polynomial} and determining the maximum eigenvalue, we derive the following expression for the global reproduction number of two orchards with ACP dispersal and control parameter:

\begin{equation}\label{eq:Ro_global}
R_{g0}  = \rho(K)  =  \dfrac{1}{\sqrt{2}}\sqrt{ b + \sqrt{b^{2} + 4c}} \, ; 
\end{equation}
with $b = (R^{2}_{10} + (1+p_{12}\delta))\cdot R^{2}_{20} $ and $c = R^{2}_{10}\cdot R^{2}_{20}$. This allows us to observe the explicit dependence on the local basic reproduction number of orchard 1 and 2, as well as the implicit dependence on the controlling parameters denoted by the hat and bold typography.

In order to evaluate the impact of sensitivity analysis on the global and local reproduction numbers (Eq. \eqref{eq:Ro_global}) in response to variations in control parameters, we perform a global sensitivity analysis using Latin Hypercube Sampling (LHS) and partial rank correlation coefficients (PRCCs) methodologies.
\begin{figure}[t]
\begin{center}
\includegraphics[width=12.5cm,height=8cm]{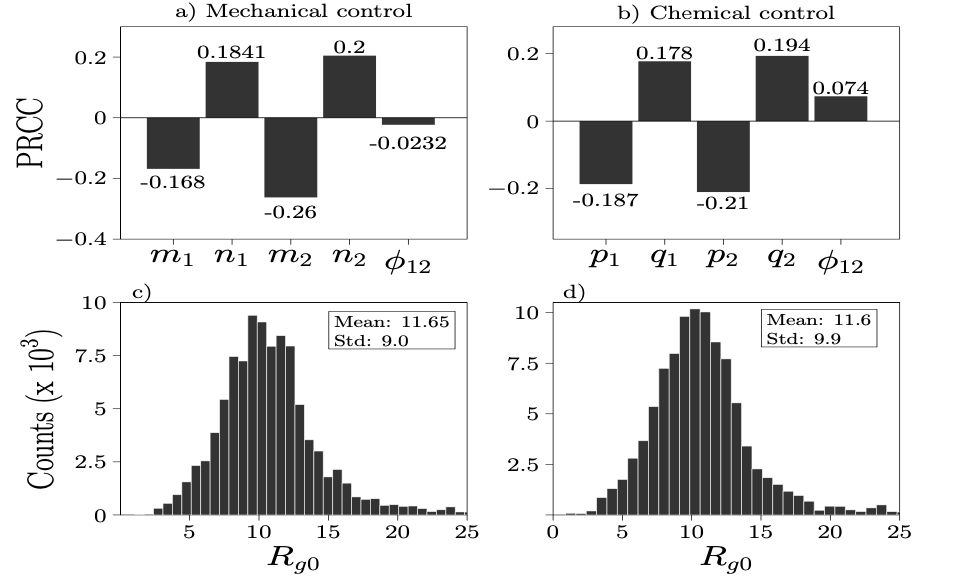}  
\caption{PRCCs for the two-orcard SAIR-SI epidemic model with mechanical and chemical control (Figures a and b). Histograms of $R_{g0}$ obtained from LHS using a sample size of $10^{5}$ for parameter values for a mechanical and chemical control strategy.}
\label{fig:PRCC_Mechanical_Ctr}
\end{center}
\end{figure}

LHS is a stratified sampling technique that enables efficient and comprehensive exploration of parameter space \cite{Blower_1994, Seaholm_1988, Helton_2005}. This sampling approach ensures equal sampling of each parameter value across its range, minimizing the risk of overlooking important parameter combinations.

Instead, PRCCs are valuable in determining the most influential control parameters and quantifying the impact of changes in each parameter on the output of the system. This is achieved through sampling using the LHS method \cite{Bidah_2020}. PRCCs allow for the assessment of the linear relationship between the global and local reproduction numbers and each input control parameter, while controlling for the effects of other parameters \cite{Marino_2008}. A PRCC value of zero indicates no linear correlation between the model output and the input parameter, while a value of -1 or 1 signifies a strong negative or positive correlation, respectively. 

To obtain uniformly distributed values for control parameters according to the LHS method, we generated $10^{5}$ samples using the \emph{pyDOE} Python library \cite{pyDOE}. For each sample (\emph{i.e.}, combination of parameter values), we evaluated the global reproduction number (Eq. \eqref{eq:Ro_global}) and recorded the resulting outputs.

We utilized the LHS method to compute the partial rank correlation coefficients PRCCs of the model outputs with the aim to perform the sensitivity analyses. For a mechanical control strategy, we set $p_{i} = q_{i} = 0$ (for $i = 1,2$,  \emph{i.e}., chemical control parameters are set to zero) and assessed the impact on the global reproduction numbers of parameters $m_{i}$ , $n_{i}$ and ACP dispersal $p_{12}$. The results are shown in Figure \ref{fig:PRCC_Mechanical_Ctr}-(a), revealing that the first control parameter $m_{1}$, associated with reducing the probing rate of tress $b_{i}$, exhibits a negative correlation with the global reproduction numbers. Conversely, the control parameter $n_{i}$, which is associated with orchard vigilance, shows a positive correlation, meaning that an increase in this parameter implies a reduction in orchard vigilance and an increase in HLB transmission. 

On the other hand, in Figure \ref{fig:PRCC_Mechanical_Ctr}-(b) we show the PRCC for a chemical control strategy, where we set $p_{i} = q_{i} = 0$ (for $i= 1,2$, \emph{i.e}., mechanical control parameters are set to zero). We observe a negative correlation between the global reproduction number and the control parameter $p_{i}$, associated with the maximum abundance of ACP vector per citrus trees. And, a positive correlation with the control parameter $q_{i}$, associated with the increase of ACP mortality. Figures \ref{fig:PRCC_Mechanical_Ctr}-(c) and \ref{fig:PRCC_Mechanical_Ctr}-(d) display the global reproduction number distribution for mechanical and chemical control strategies, respectively. 


\section{Cost and effectiveness of the control strategies.}

In this section, we present a mathematical framework to evaluate the cost associated with intervention strategies in an orchard, specifically focusing on mechanical and chemical control strategies. Our aim is to devise a numerical approach for estimating the cost of these strategies and their effectiveness, taking into account the mobility parameter $\phi_{12}$ and the local reproduction number $R_{i0}$ (Eq. \eqref{eq:Rio}). By employing this framework, we can assess the financial implications of different control measures and their impact on mitigating the spread of the HLB disease in the orchard. 

\subsection{Cost.}

Based on \cite{FZhang_2021}, we define the mechanical control variable $\mathbf{X}$, which is associated with the control parameters $m_{1}$ and $n_{1}$ for orchard 1, and $m_{2}$ and $n_{2}$ for orchard 2. Similarly, we have the chemical control variable $\mathbf{Y}$, which is associated with the chemical control parameters $p_{1}$ and $q_{1}$ for orchard 1, and $p_{2}$ and $q_{2}$ for orchard 2. In this context, we define the cost of using a mechanical control strategy as follows:
\begin{equation}\label{eq:cost_mechanical_ctrl}
\begin{array}{lll}
C_{\mathbf{X}} &=& \displaystyle \int_0^\infty \Big( M_{1}(A_{1\tau}(t) + I_{1\tau}(t)) + M_{2}(A_{2\tau}(t) + I_{2\tau}(t)) \Big) dt + \\ 
\\
&+& \displaystyle \int_{T_{i}}^{T_{f}} \Big( x_{1} m^{2}_{1}(t) + x_{2} n^{2}_{1}(t) + x_{3} m^{2}_{2}(t) + x_{4} n^{2}_{2}(t)\Big)dt \, ;
\end{array}
\end{equation}

\noindent and, for the chemical control strategy, we define its cost as:
\begin{equation}\label{eq:cost_chemical_ctrl}
\begin{array}{lll}
C_{\mathbf{Y}} &=& \displaystyle \int_0^\infty \Big( C_{1}(A_{1\tau}(t) + I_{1\tau}(t)) + C_{2}(A_{2\tau}(t) + I_{2\tau}(t)) \Big) dt + \\ 
\\
&+& \displaystyle \int_{T_{i}}^{T_{f}} \Big( y_{1} p^{2}_{1}(t) + y_{2} q^{2}_{1}(t) + y_{3} p^{2}_{2}(t) + y_{4} q^{2}_{2}(t)\Big)dt \, .
\end{array}
\end{equation}

The time interval during which these control strategies are applied is denoted by $T_{i}$ and $T_{f}$, representing the initial and final time, respectively.

The coefficients $M_{i}$ and $C_{i}$ ($i=1,2$) are positive constants, ranging from zero to one, that represent the weights assigned to the number of infectious asymptomatic and symptomatic citrus trees in orchard 1 and 2, respectively. The coefficients $x_{i} \in (0, 1]$ ($i=1,2,3,4$) are weight constants that denote the costs associated with the mechanical control strategy, which includes the use of protective barriers, yellow-colored traps, increased orchard vigilance, and the effectiveness of identifying infectious symptomatic trees. Similarly, the coefficients $y_{i} \in (0, 1]$ ($i=1,2,3,4$) represent the corresponding weight constants associated with the cost of the chemical control strategy. This strategy involves the use of external means such as fumigation, biological control, or cultural practices such as the use of alternate host plants. We assume that the costs are proportional to a quadratic form of the corresponding control functions.

In this paper, we consider all cost weights to be equal to 1, and they will remain constant throughout the subsequent analysis. It is essential to emphasize that the cost weights used in the simulations are purely for theoretical purposes and serve to illustrate the proposed control strategies in this article. 

A relevant aspect to emphasize regarding our definition of costs (Eqs. \eqref{eq:cost_mechanical_ctrl} and \eqref{eq:cost_chemical_ctrl}) is the integration range in the first integral, which corresponds to the cumulative fraction of infectious asymptomatic and symptomatic citrus trees. This integration is performed from zero to infinity. The rationale behind this choice is to capture the cumulative effect of all infections over the entire duration of the epidemic. By integrating from zero to infinity in the first integral, we are considering the long-term impact of infections, encompassing both early and late stages of the epidemic. In contrast, the second integral in the cost function is limited to the time period during which the control strategy is applied. This reflects the fact that the cost calculation takes into account the specific time period during which the control strategies are actively employed, ensuring that the cost calculation is focused on the specific intervention period.

In \ref{appendix_A}, we derive an explicit expression for the integral of the compartments $A_{i}(t) + I_{i}(t)$ from zero to infinity, utilizing the model equations (Eqs. \EqsModel). This expression allows us to rewrite the cost function in a more concise form as:
\begin{adjustwidth}{-0.5cm}{-0.5cm} 
\begin{align*}
C_{\mathbf{X}} = \Gamma_{1} R_{1\tau}(\infty) + \Gamma_{2} R_{2\tau}(\infty) + \int_{T_{i}}^{T_{f}} \Big( x_{1} m^{2}_{1}(t) + x_{2} n^{2}_{1}(t) + x_{3} m^{2}_{2}(t) + x_{4} n^{2}_{2}(t)\Big) dt , \\
\\
C_{\mathbf{Y}} = \Gamma_{1} R_{1\tau}(\infty) + \Gamma_{2} R_{2\tau}(\infty) + \int_{T_{i}}^{T_{f}} \Big( y_{1} p^{2}_{1}(t) + x_{2} q^{2}_{1}(t) + y_{3} p^{2}_{2}(t) + y_{4} q^{2}_{2}(t)\Big) dt ;
\end{align*}
\end{adjustwidth}

\noindent where  $\Gamma_{i} = (\sigma_{i} + \bmh{r}_{i\tau})/(\bmh{r}_{i\tau}\cdot (\sigma_{1} + \mu_{1\tau}))$ (with $i=1,2$), and we make the assumption that there are no infectious symptomatic trees at the beginning of the disease, i.e., $I_{i\tau}(0)=0$ in both orchards. Here, $R_{i}(\infty) = \lim_{t \rightarrow \infty} R_{i}(t)$ indicates the final count of rogue trees throughout the entire epidemic duration.

\subsection{Effectiveness.}
We define the effectiveness of the control strategy $X$ (or Y) in terms of the final size of the disease, \textit{i.e} $S_{i}(\infty) = \lim_{t \rightarrow \infty} S_{i}(t)$ , which corresponds to the total number of infected trees at the end of the epidemic. This measure allows us to evaluate the overall impact and success of the control strategy in mitigating the spread of HLB. By focusing on the final size, we capture the cumulative effect of the control strategies implemented over the entire course of the epidemic. In \ref{appendix_B}, we show that the explicit expression for the orchard $1$ final size, and therefore the effectiveness in orchard 1 (Eq. \eqref{Final_Size_Trees_Orchard_1}), is given by:
\begin{flalign}\label{Effectivenes_Orchard_1}
Ef_{1} & = S_{1\tau}(0)\exp\left( -\dfrac{R^{2}_{10}\phi^{2}_{11}}{N_{1\tau}}\cdot R_{1\tau}(\infty) - \theta_{12} \dfrac{R^{2}_{20}\bmh{\beta}_{11}}{N_{1\tau}\bmh{\beta}_{21}} \cdot \phi_{11}\phi_{12} \cdot R_{2\tau}(\infty) \right) 
\, ;
\end{flalign}

\noindent and, for the second orchard (Eq. \eqref{Final_Size_Trees_Orchard_2}), we obtain:
\begin{flalign}\label{Effectivenes_Orchard_2}
Ef_{2} & = S_{2\tau}(0)\exp\left(-\dfrac{R^{2}_{20}}{N_{2\tau}} ( 1 + \theta_{12} p^{2}_{12}) \cdot R_{2\tau}(\infty)   
 -\dfrac{R^{2}_{10}\bmh{\beta}_{21}}{N_{2\tau}\bmh{\beta}_{11}}\phi_{11}\phi_{12} \cdot R_{1\tau}(\infty)  \right)  \, ;
\end{flalign}

\noindent where we are setting $I_{i} = 0$ for i=1,2;  $\theta_{12} = N^{0}_{1v}\bmh{\mu}_{2v}/ N^{0}_{2v}\bmh{\mu}_{1v}$ and $R_{i0}$ is the local basic reproduction number given in Eq. \eqref{eq:Rio}.

\section{Control strategies based on genetic algorithms.}

A genetic algorithm is a computational optimization technique that is inspired by the process of natural selection and evolution \cite{Holland_1992, Mitchell_1998}. It involves the use of genetic operators, such as mutation and crossover, to evolve a population of candidate solutions towards an optimal solution. Through successive generations of selection, reproduction, and mutation, genetic algorithms are capable of identifying high-quality solutions to complex optimization problems. By mimicking the principles of evolution, genetic algorithms effectively explore the solution space and converge towards optimal solutions, making them a powerful tool for solving intricate optimization challenges.

In this section, we utilized a genetic algorithm for determining the optimal combination of control parameters to minimize the number of rogued citrus trees and reduce the overall cost in the context of HLB disease. Our approach utilizes the system of differential equations that describe the dynamics of the two-orchard model, as well as the cost functions defined in Eqs. \eqref{eq:cost_mechanical_ctrl} and \eqref{eq:cost_chemical_ctrl}. By integrating these components, we aim to identify the most effective control strategy that achieves the desired outcome of minimizing tree roguing while considering the associated costs.


\subsection{Control objective.}

The control objective for the genetic algorithm is to systematically determine the optimal control parameters that allow the system of Eqs. \EqsModel to evolve in a way that, at the conclusion of the HLB epidemic, the total number of roguing trees $R_{1\tau}(\infty) + R_{2\tau}(\infty)$ achieves a desired total number $R^{obj}$.

To fulfill this objective, the genetic algorithm aims to minimize the following objective function:

\begin{equation}\label{eq:objetive_function}
J_{X} = w_{1} \cdot C_{X} + w_{2} \cdot (Ef_{1} + Ef_{2}) \, ,
\end{equation}

where $J_{X}$ represents the objective function for the mechanical control strategy, $C_{X}$ is the cost of implementing this control strategy, and $Ef_{1}$ and $Ef_{2}$ denote the efficiencies of the control measures in each orchard. The weights $w_{1}$ and $w_{2}$ are used to adjust the relative importance of the cost and control efficiency terms in the objective function. 


\subsection{Algorithm description.}

In this section, we outline the genetic algorithm by utilizing the mechanical control strategy $X$ as an illustrative example. Specifically, we aim to determine the optimal values for the control parameters $m_{1}$, $n_{1}$, $m_{2}$, $n_{2}$, as well as the corresponding cost function $C_{X}$ defined in Eq. \eqref{eq:cost_mechanical_ctrl}. It is important to note that the same principles can be extended to the chemical control strategy by substituting the control parameters with $p_{1}$, $q_{1}$, $p_{2}$, $q_{2}$, and the cost function with $C_{Y}$ defined in Eq. \eqref{eq:cost_chemical_ctrl}. With this in mind, the genetic algorithm can be outlined through the following steps:

\textbf{Step 0}: The non-controlling parameters are initialized with their baseline values as specified in Table \eqref{table_parameters}. The initial conditions for both orchards are determined using Eq. \eqref{eq:init_cond_orchar_1} for orchard 1 and Eq. \eqref{eq:init_cond_orchar_2} for orchard 2. To initiate the genetic algorithm process, we generate a set of $M$ combinations of control parameter values $m_{1}$, $n_{1}$, $m_{2}$, $n_{2}$ uniformly distributed between zero and one. Each combination of parameters is referred to as an individual $k \in \{1,2,\ldots, M \}$ of the generation. For each individual, we evaluate the two-orchard model described by Eqs. \EqsModel and record the values of $R_{i\tau}(\infty)$ and $S_{i\tau}(\infty)$. These values are then used to calculate the cost function $C_{X}$ (Eq. \eqref{eq:cost_chemical_ctrl}), control effectiveness $Ef_{i}$, and the objective function $J$ (Eq. \eqref{eq:objetive_function}).

\textbf{Step 1}: (\emph{Best fitness}) In each successive generation, we select from the entire population, a single individual $k^{ftr}$ with the best \comillas{fitness}, that is, we choose the parameter combinations from the previous generation that satisfy the condition $|(R_{1\tau}(\infty) + R_{2\tau}(\infty)) - R^{obj}| \approx 0$, with a preference for those with the minimum value of the objective function $J$. The combination of parameters that constitute the individual $k^{ftr}$ are labeled as $m_{1}^{ftr}$, $n_{1}^{ftr}$, $m_{2}^{ftr}$, $n_{2}^{ftr}$. Starting from this individual, we generate the remaining population of parameter combinations for the current generation, following the procedure outlined below:

\textbf{Step 2}: (\emph{Segment crossover}) From the population of $M-1$ individuals (excluding the chosen individual $k^{ftr}$ from Step 1), we select $N$ individuals with the highest fitness, named as $k^{mth}_{i}$ for $i=1, 2, \ldots, N$. We denote as $m_{1i}^{mth}$, $n_{1i}^{mth}$, $m_{2i}^{mth}$, $n_{2i}^{mth}$, the combination of parameters associated to the individual $k^{mth}_{i}$. For each selected individual $k^{mth}_{i}$, we perform crossover by combining their genetic information with the best selection in Step 1 $k^{ftr}$. In specific, two random points are generated, and the segments between these points in the father chromosomes ($m_{1}^{ftr}$ for example) are swapped with the mother individual ($m_{1i}^{mth}$ for example) and create the first child; next, the second child is created by cutting off the segments of the mother and pasting them into the father. 

\textbf{Step 3}: (\emph{Mutation and elitism}) In this step of the genetic algorithm, we introduce mutation to further diversify the offspring population. Randomly, selected segments (genes) from the descendant individuals undergo mutation, resulting in new values for these segments. The probability of mutation is set $0.1$, ensuring a balance between exploration and exploitation in the search process. The mutated values are recorded as the parameters for the offspring individuals. Additionally, the elitism process is applied to preserve the best solutions obtained so far. the individual $k^{ftr}$ is retained without any modification, ensuring that the top-performing individuals from the previous generation are included in the current generation. 

\textbf{Step 4}: (\emph{Survival}) This step is used to select the remaining individuals $k_{j}^{surv}$ (for $j=1,....,S$) for the mating pool after the elite individuals have been selected. Individuals with lower fitness scores, who were not selected for elitism, can be included in the mating pool through this selection process. This process increases the diversity of the population and allows for exploration of the search space beyond local optima. The selection process begins by randomly selecting pairs of individuals from the remaining pool and performing the segment crossover operation described in \textbf{Step 2}. This process continues until a new generation of $M$ individuals is generated.

\textbf{Step 5}: After generating all combinations of parameters of the new generation, we proceed to evaluate each combination using the two-orchard model described by Eqs. \EqsModel , \emph{i.e.}, for each combination, we record the values of $R_{i\tau}(\infty)$ and $S_{i\tau}(\infty)$, which are crucial for calculating the cost function $C_{X}$ (Eq. \eqref{eq:cost_chemical_ctrl}), control effectiveness $Ef_{i}$, and the objective function $J$ (Eq. \eqref{eq:objetive_function}). 

In the subsequent section, we utilize the genetic algorithm, outlined previously, to optimize the control parameters, aiming to minimize the objective function $J$ (Eq. \eqref{eq:objetive_function}) and achieve the desired number of roguing trees ($R^{obj}$).

\subsubsection{Numerical results.}

\begin{figure}[ht]
\begin{center}
\includegraphics[width=13cm,height=9.5cm]{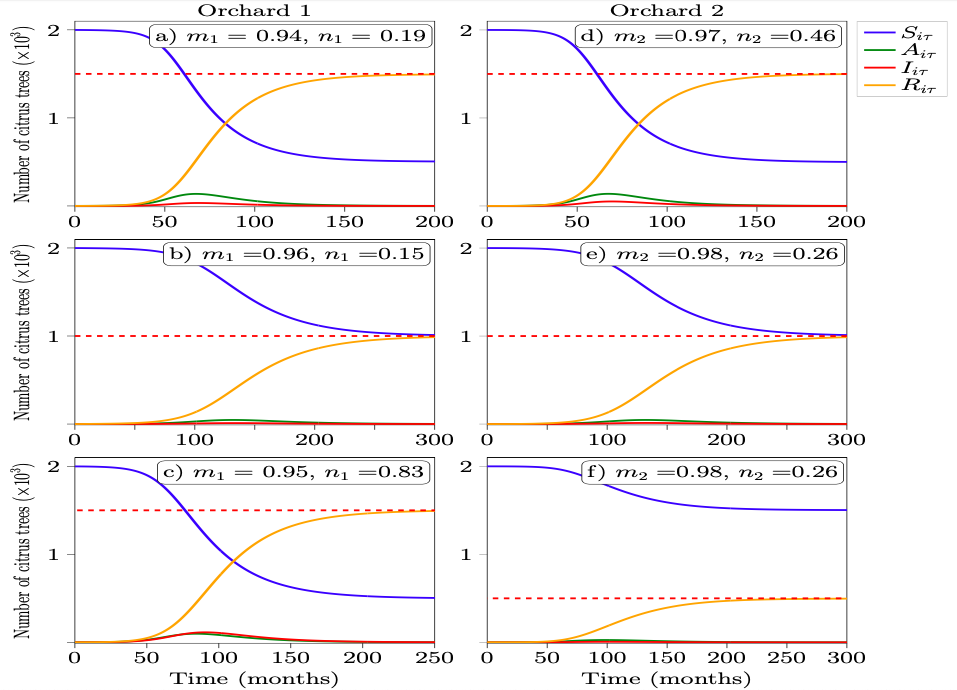}  
\caption{Numerical solution to the two-orchard SAIR-SI model (Eqs. \EqsModel) for three examples using different values of the control parameters $m_{i}$ and $n_{i}$ (with $i=1,2$)  obtained through a genetic algorithm for a mechanical control strategy. The red dashed lines represent the values of $R_{1\tau}$, and $R_{2\tau}$, confirming the effectiveness of the genetic algorithm in achieving the desired objectives for the mechanical control strategy.}
\label{Fig:mechanical_strategy_Sol}
\end{center}
\end{figure}

In order to evaluate the functionality of the algorithm, our objective was to determine the optimal combination of control parameters that would result in a final size of the SAIR-SI model equal to $R^{obj} \rightarrow 0.75 N_{1\tau}$, $R^{obj} \rightarrow 0.5 N_{1\tau}$, or $R^{obj} \rightarrow 0.25 N_{1\tau}$, depending on the desired percentage of roguing trees at the end of the HLB epidemic (at what follows, we assume that $N_{1\tau} = N_{2\tau} = 2,000$). The main goal was to achieve a specific total number of trees to be removed from both orchards, while minimizing the objective function \eqref{eq:objetive_function}. Additionally, we introduced an extra constraint to the algorithm, which involved selecting a specific number of roguing trees for each individual orchard. For instance, if each orchard contained $2,000$ trees and $R^{obj}=3,000$, we could enforce the constraint that $R_{1\tau}(\infty) = 2,000$ trees should be rogued from the first orchard and $R_{2\tau}(\infty) = 1,000$ trees from the second orchard.

In this sense, we propose the following three examples for a mechanical and chemical control strategy:

\begin{itemize}
    \item[$\ast$]  Example 1: $R^{obj} = 0.75 \cdot N_{1\tau} = 3,000$ constrained to $R_{1\tau}(\infty) = 1,500$, $R_{2\tau}(\infty) = 1,500$ and $min\{J_{X}\}$.
    \item[$\ast$]  Example 2: $R^{obj} = 0.5 \cdot N_{1\tau}  = 2,000$ constrained to $R_{1\tau}(\infty) = 1,000$, $R_{2\tau}(\infty) = 1,000$ and $min\{J_{X}\}$.
    \item[$\ast$]  Example 3A: $R^{obj} = 0.5 \cdot N_{1\tau} = 2,000$ constrained to $R_{1\tau}(\infty) = 1,500$, $R_{2\tau}(\infty) = 500$ and $min\{J_{X}\}$. This example is explored for a mechanical control strategy. 
    \item[$\ast$]  Example 3B: $R^{obj} = 0.25 \cdot N_{1\tau} = 1,000$ constrained to $R_{1\tau}(\infty) = 500$, $R_{2\tau}(\infty) = 500$ and $min\{J_{Y}\}$. This example is explored for a chemical control strategy. 
\end{itemize}

In all the previously mentioned scenarios, we made the assumption that the ACP emigration parameter, denoted as $\phi_{12}$, is fixed at a value of $0.35$. This means that only $35$ \% of the ACP population in orchard 1 migrates to orchard 2. It is important to note that the chosen value for $\phi_{12}$ is a simplification and may not accurately represent the real-world ACP migration rates in every specific situation. However, it serves as a useful approximation to investigate the effects of the other control parameters on the overall dynamics of the system.

\begin{figure}[ht]
\begin{center}
\includegraphics[width=13cm,height=9.5cm]{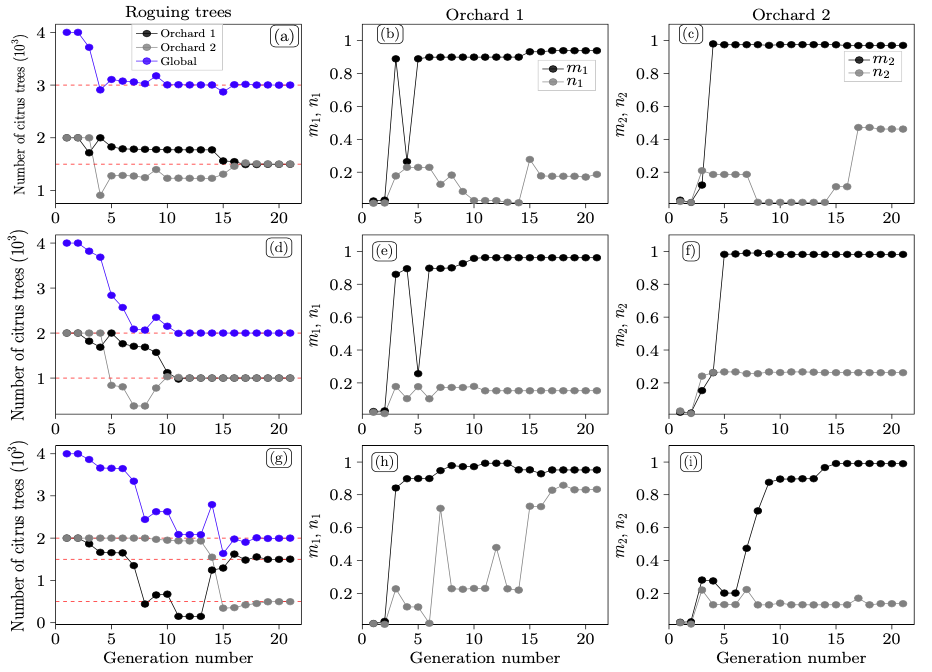}  
\caption{The progressive changes in the parameter values $m_i$ and $n_i$ over the generations during the genetic algorithm optimization process for a mechanical control strategy. (a) the evolution of the control objectives $R^{obj}$, $R_{1\tau}$, and $R_{2\tau}$ (represented with the red dashed lines) for the best combination of parameters in each generation. (b) and (c) depict the evolution of $m_{1}$ and $n_{1}$ for orchard 1, and $m_{2}$ and $n_{2}$ for orchard 2, respectively.} 
\label{Fig:mechanical_strategy_Genetic_Algorithm}
\end{center}
\end{figure}

In Figure \eqref{Fig:mechanical_strategy_Sol} we observe the solution to the two-orchard model with the parameter values $m_{i}$ and $n_{i}$ (with $i=1,2$) finding with the genetic algorithm for each one of the aforementioned examples. We observe that the parameter $m_{i}$ requires high values, resulting in a cost of $C_{X} = 19,388$ for example 1, $C_{X} = 12,397$ for example 2, and $C_{X} = 19,311$ for example 3, highlighting the significant cost associated with implementing this mechanical strategy.

Furthermore, Figure \eqref{Fig:mechanical_strategy_Genetic_Algorithm} illustrates the progressive changes in the parameter values $m_i$ and $n_i$ as the number of generations increases during the genetic algorithm optimization process. Additionally, the figure displays (in red dashed lines) the corresponding control objectives $R^{obj}$, $R_{1\tau}$, and $R_{2\tau}$ for the best parameter combination in each generation. We observe that the genetic algorithm requires 15 generations to achieve the objective control for the first and third examples, while it only takes 10 generations for the second example.
\begin{figure}[ht]
\begin{center}
\includegraphics[width=13cm,height=9.5cm]{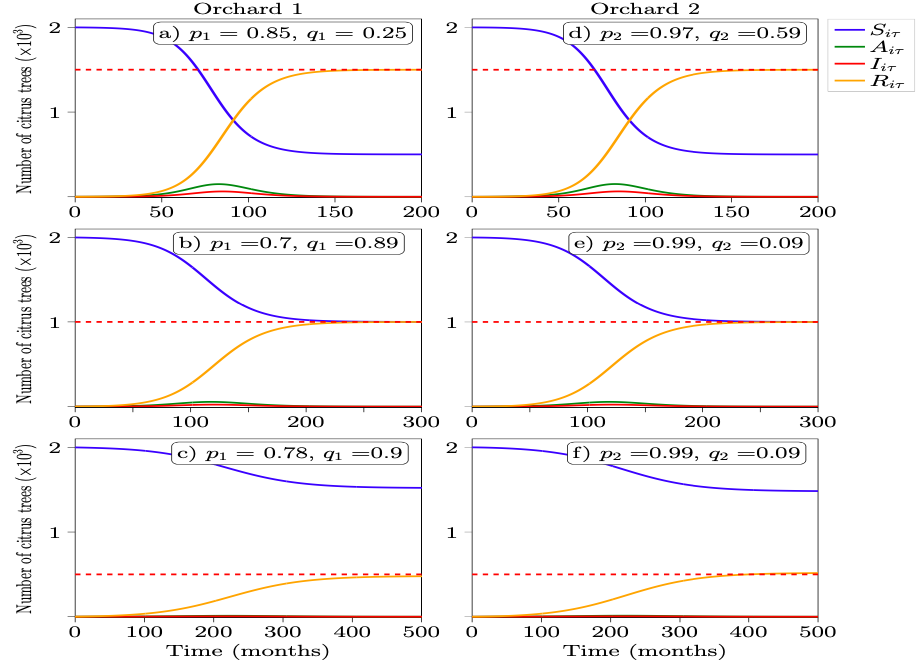}  
\caption{Numerical solution to the two-orchard SAIR-SI model (Eqs. \EqsModel) for three examples using different values of the control parameters $p_{i}$ and $q_{i}$ (with $i=1,2$), obtained through a genetic algorithm for a mechanical control strategy. The red dashed lines represent the values of $R_{1\tau}$, and $R_{2\tau}$, confirming the effectiveness of the genetic algorithm in achieving the desired objectives for the chemical control strategy.}
\label{Fig:chemical_strategy_Sol}
\end{center}
\end{figure}

In Figure \eqref{Fig:chemical_strategy_Sol}, we present the numerical solution to the two-orchard model using the parameter values $p_{i}$ and $q_{i}$ (where $i=1,2$) obtained through the genetic algorithm. The associated cost for implementing this chemical control strategy is $C_{Y} = 21,106$ for the first example, $C_{Y} = 14,093$ for the second example, and $C_{Y} = 7,109$ for the third example.

\begin{figure}[ht]
\begin{center}
\includegraphics[width=13cm,height=9.5cm]{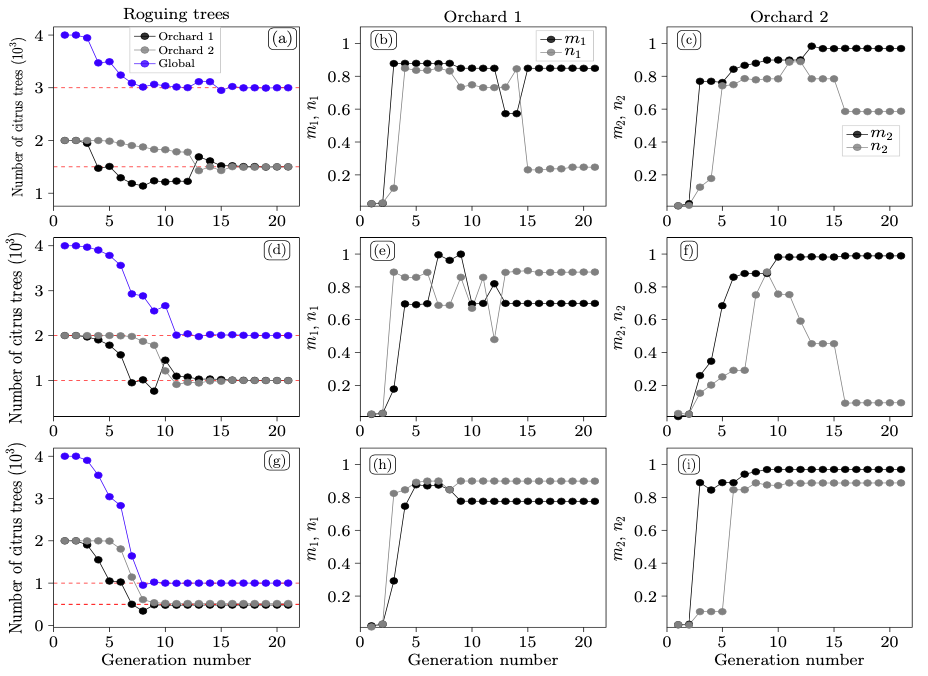}  
\caption{The progressive changes in the parameter values $p_{i}$ and $q_{i}$ over the generations during the genetic algorithm optimization process for a chemical control strategy. (a) the evolution of the control objectives $R^{obj}$, $R_{1\tau}$, and $R_{2\tau}$ (represented with the red dashed lines) for the best combination of parameters in each generation. (b) and (c) depict the evolution of $p_{1}$ and $q_{1}$ for orchard 1, and $p_{2}$ and $q_{2}$ for orchard 2, respectively.} 
\label{Fig:chemical_strategy_Genetic_Algorithm}
\end{center}
\end{figure}

The progressive changes in the parameter values $p_{i}$ and $q_{i}$ (with $i=1,2$) with respect to the generation number is shown in Figure \eqref{Fig:chemical_strategy_Genetic_Algorithm}. For the first example, the genetic algorithm requires 15 generations to achieve the objective control, while for the second and third examples, it takes at most 10 generations.

\section{Conclusions.}

This study demonstrates the effectiveness of utilizing genetic algorithms to optimize control parameters in mitigating the spread of Huanglongbing (HLB) in citrus orchards. The developed two-orchard model, considering the dispersal of the Asian Citrus Psyllid (ACP), provides valuable insights into the efficacy of mechanical and chemical control strategies. By incorporating ACP mobility, the model captures the dynamics of HLB more realistically and allows for sensitivity analysis of model parameters.

The results obtained through the genetic algorithms showcase the successful identification of optimal control parameters for each strategy. The derived mathematical expression for the global reproduction number ($R_{0}$) enables a comprehensive assessment of the model sensitivity to ACP mobility. Furthermore, the cost functions and efficiency analysis provide a quantitative evaluation of the control strategies based on the final size and individual $R_{0}$ of each orchard.

This research emphasizes the significance of optimizing control parameters in disease management within the agricultural sector. The use of genetic algorithms not only improves our understanding of HLB control but also serves as a foundation for future studies in developing disease control strategies based on this innovative approach. These findings contribute to informed decision-making and the implementation of effective control measures against HLB in citrus orchards, ultimately aiding in the preservation and sustainability of citrus production.


\appendix

\section{The cumulative fraction of ACP vectors.}\label{appendix_A}

Our objective is to obtain a mathematical expression for the integral of compartments $A_{i\tau}(t) + I_{i\tau}(t)$ from zero to infinity using Eqs. \EqsModel. This integral represents the cumulative fraction of citrus trees that have been infected by HLB bacteria throughout the entire duration of the epidemic.

First, we can obtain an expression for the integral of compartment $I_{1\tau}(t)$ by integrating both sides of the equation for the roguing compartment $R_{1\tau}$ in Eq. \eqref{eq:orchard1_EqR1t} from 0 to infinity, considering its initial condition $R_{\tau}(0) = 0$ (which indicates that no trees were removed at the beginning of the disease):
\begin{equation}\label{eq:Int_I1t}
\int_0^\infty I_{1\tau}(t) dt = \frac{R_{1\tau}(\infty)}{\bmh{r}_{1\tau}} - \frac{\mu_{1\tau}}{\bmh{r}_{1\tau}}\int_0^\infty A_{1\tau}(t)dt \, .
\end{equation}

Similarly, we can obtain an expression for the integral of compartment $I_{2\tau}(t)$ by integrating both sides of the equation for the roguing compartment $R_{2\tau}$ in Eq. \eqref{eq:orchard2_EqR2t}:
\begin{equation}\label{eq:Int_I2t}
\int_0^\infty I_{2\tau}(t) dt = \frac{R_{2\tau}(\infty)}{\bmh{r}_{2\tau}} - \frac{\mu_{2\tau}}{\bmh{r}_{2\tau}}\int_0^\infty A_{2\tau}(t)dt \, ;
\end{equation}

\noindent where $R_{i\tau}(\infty) = lim_{t \rightarrow \infty} R_{i\tau}$, for $i=1,2$ is called the final size for roguing trees. 

On the other hand,  by summing the compartments   $I_{1\tau}$ and $R_{1\tau}$ (Eqs.   \eqref{eq:orchard1_EqI1t} and \eqref{eq:orchard1_EqR1t}), we get: 
\[
\dfrac{dR_{1\tau}}{dt} + \dfrac{dI_{1\tau}}{dt} = (\sigma_{1} + \mu_{1\tau})A_{1\tau} \, ;
\]
\noindent from which we get, after integrating both sides and considering that $R_{1\tau}(0) = 0$:
\begin{equation}\label{eq:Int_A1t}
    \int_0^\infty  A_{1\tau}(t) dt  = \dfrac{R_{1\tau}(\infty)}{\sigma_{1}+\mu_{1\tau}} - \dfrac{I_{1\tau}(0)}{\sigma_{1}+\mu_{1\tau}} \, ;
\end{equation}

\noindent and doing the same for compartments $I_{2\tau}$ and $R_{2\tau}$ (from Eqs. \eqref{eq:orchard2_EqI2t} and \eqref{eq:orchard2_EqR2t}), we get:
\begin{equation}\label{eq:Int_A2t}
    \int_0^\infty  A_{2\tau}(t) dt  = \dfrac{R_{2\tau}(\infty)}{\sigma_{2}+\mu_{2\tau}} - \dfrac{I_{2\tau}(0)}{\sigma_{2}+\mu_{2\tau}} \,.
\end{equation}

Then, from Eq. \eqref{eq:Int_I1t} and Eq. \eqref{eq:Int_A1t}  we get:
\begin{equation}\label{eq:sum_in_A1_I1}
    \begin{array}{lcl}
       \int_0^\infty ( A_{1\tau}(t) + I_{1\tau}(t) ) dt   & =& \dfrac{R_{1\tau}(\infty)}{\bmh{r}_{1\tau}} + \Big( 1 - \dfrac{\mu_{1\tau}}{\bmh{r}_{1\tau}}\Big) \int_0^\infty  A_{1\tau}(t) dt \\
          \\
          & = & \dfrac{R_{1\tau}(\infty)}{\bmh{r}_{1\tau}}\dfrac{(\sigma_{1} + \bmh{r}_{1\tau})}{(\sigma_{1} + \mu_{1\tau})} - \dfrac{(\bmh{r}_{1\tau} - \mu_{1\tau} )}{\bmh{r}_{1\tau}(\sigma_{1} + \mu_{1\tau})} I_{1\tau}(0).
    \end{array}
\end{equation}

Similarly, by using Eq. \eqref{eq:Int_I2t} and Eq. \eqref{eq:Int_A2t}, we can obtain the same expression for orchard 2 by replacing the orchard label 1 with 2 ($1 \rightarrow 2$) in the above equation.

\section{Final size for citrus trees.}\label{appendix_B}

Consider the equation for the compartment $S_{1\tau}$ in Eq. \eqref{eq:orchard1_EqS1t} for orchard $1$. It could be rewritten as 
\[
\int_0^\infty \dfrac{dS_{1\tau}}{S_{1\tau}} = -\dfrac{\bmh{\beta}_{11}\cdot \phi_{11}}{N_{1\tau}}\int_0^\infty I_{1v}(t)dt \, ,
\]

From which we obtain:
\begin{equation}\label{eq:Log_S1t_S10}
ln\Biggl( \dfrac{S_{1\tau}(\infty)}{S_{1\tau}(0)}\Biggr) =-\dfrac{\bmh{\beta}_{11} \phi_{11}}{N_{1\tau}} \int_0^\infty I_{1v}(t)dt 
\end{equation}

\noindent where $S_{1\tau}(\infty) = \lim_{t \rightarrow \infty} S_{1\tau}(t)$  and $S_{1\tau}(0)$ is the total number of susceptible citrus trees at the beginning of the HLB epidemic in the orchard 1. Then:
\begin{equation}\label{eq:S1_Infty}
    S_{1\tau}(t) = S_{1\tau}(0)\cdot \exp\left(-\dfrac{\bmh{\beta}_{11} \phi_{11}}{N_{1\tau}} \int_0^\infty I_{1v}(t)dt \right)
\end{equation}

The objective now is to derive an expression for the integral of the compartment $I_{iv}(t)$ from zero to infinity. Since $S_{1v}(t) = N_{1v} -  I_{1v}(t)$, then, from Eq. \eqref{eq:orchard1_EqI1v} we get:
\[
\dfrac{dI_{1v}}{dt} = \bmh{\beta}_{1v} \phi_{11} (N_{1v} - I_{1v})\dfrac{A_{1\tau} + I_{1\tau}}{N_{1\tau}} + \bmh{\beta}_{2v} \phi_{12}(N_{1v} -  I_{1v})\dfrac{A_{2\tau} + I_{2\tau}}{N_{2\tau}}  - \bmh{\mu}_{1v}I_{1v} \,.
\]

\noindent Solving for $A_{1\tau} + I_{1\tau}$ we get: 

\begin{equation}\label{eq:At_plus_It}
\begin{array}{lll}
A_{1\tau} + I_{1\tau} &=& \dfrac{N_{1\tau}}{\bmh{\beta}_{1v}\phi_{11}N_{1v}(1-I_{1v}/N_{1v})}\cdot \dfrac{dI_{1v}}{dt} - \dfrac{N_{1\tau}}{\bmh{\beta}_{1v}\phi_{11}} \dfrac{\bmh{\beta}_{2v}\phi_{12}}{N_{2\tau}}\cdot(A_{2\tau} + I_{2\tau}) + \\
\\
& & \dfrac{N_{1\tau}}{\bmh{\beta}_{1v}\phi_{11}N_{1v}(1-I_{1v}/N_{1v})}\cdot \bmh{\mu}_{1v}I_{1v} \, , \\ 
\\
&\approx&  \dfrac{N_{1\tau}}{\bmh{\beta}_{1v}\phi_{11}N_{1v}}\cdot \dfrac{dI_{1v}}{dt} - \dfrac{N_{1\tau}}{N_{2\tau}} \dfrac{\bmh{\beta}_{2v}\phi_{12}}{\bmh{\beta}_{1v}\phi_{11}}\cdot(A_{2\tau} + I_{2\tau}) + \\ 
\\
& & \dfrac{N_{1\tau}\bmh{\mu}_{1v}} {\bmh{\beta}_{1v}\phi_{11}N_{1v}}\cdot I_{1v} \, ; 
\end{array}
\end{equation}

\noindent where we approximate the term $1/(1-I_{1v}/N_{1v})$  by a series expansion and ignore higher order terms (assuming that $I_{1v}/N_{1v} < 1$).

By integrating both sides from zero to infinity and considering that $I_{1v}(t)$ satisfies the following conditions: a) it converges to zero as $t \rightarrow \infty$ ($I_{1v}(\infty) \rightarrow 0$), and b) $I_{1v}(0) = 0$, indicating that there are no symptomatic infectious ACP vectors at the beginning of the disease, we can obtain the integral expression for the derivative of $I_{1v}(t)$ as follows:
\[
\int_0^\infty \dfrac{dI_{1v}(t)}{dt}dt = I_{1v}(\infty) - I_{1v}(0) = 0 \,.
\]

Then, from Eq. \eqref{eq:At_plus_It} we obtain:
\begin{equation}\label{eq:Int_I1v}
    \int_0^\infty I_{1v}(t)dt =  \dfrac{\bmh{\beta}_{1v}\phi_{11}N_{1v}}{\bmh{\mu}_{1v} N_{1\tau}} \int_0^\infty (A_{1\tau} + I_{1\tau})dt + \dfrac{\bmh{\beta}_{2v} \phi_{12}N_{1v}}{\bmh{\mu}_{1v} N_{2\tau}} \int_0^\infty (A_{2\tau} + I_{2\tau})dt \,.
\end{equation}

Substituting Eq. \eqref{eq:sum_in_A1_I1} in  Eq. \eqref{eq:Int_I1v} and using the expression for the local basic reproduction number (Eq. \eqref{eq:Rio}) for $i=1$, we obtain:
\begin{flalign}\label{eq:Int_I1v_Final}
\int_0^\infty  I_{1v}(t) dt & \approx \dfrac{R_{10}^{2}\phi_{11}}{\bmh{\beta}_{1\tau}}\cdot \Biggl( R_{1\tau}(\infty) + \dfrac{\mu_{1\tau}-\bmh{r}_{1\tau}}{\sigma_{1}+\bmh{r}_{1\tau}}\cdot I_{1\tau}(0) \Biggr) + \nonumber \\ 
\\ \nonumber
& +  \dfrac{N^{0}_{1v}\bmh{\mu}_{2v}}{N^{0}_{2v}\bmh{\mu}_{1v}} \cdot \dfrac{R^{2}_{20}\phi_{12}}{\bmh{\beta}_{2\tau}} \Biggl( R_{2\tau}(\infty) +  \dfrac{\mu_{2\tau} - \bmh{r}_{2\tau}}{\sigma_{2} + r_{2\tau}} \cdot I_{2\tau}(0)\Biggr) \, .
\end{flalign}

Then, substituting in \eqref{eq:S1_Infty} we obtain: 
\begin{flalign}\label{Final_Size_Trees_Orchard_1}
S_{1\tau}(\infty) & = S_{1\tau}(0)\exp\left( -\dfrac{R^{2}_{10}\phi^{2}_{11}}{N_{1\tau}}\cdot \Biggl( R_{1\tau}(\infty) + \delta_{1} I_{1\tau}(0) \Biggr) \right) \times \nonumber \\ 
& \hspace{1.3cm} \times \exp\left(  - \theta_{12} \dfrac{R^{2}_{20}\bmh{\beta}_{11}}{N_{1\tau}\bmh{\beta}_{21}}
\cdot \phi_{11}\phi_{12} \cdot \Biggl( R_{2\tau}(\infty) + \delta_{2} I_{2\tau}(0)\Biggr) \right)
\, ;
\end{flalign}

\noindent where 
\[
\delta_{1} = \dfrac{\mu_{1\tau} - \bmh{r}_{1\tau}}{\sigma_{1} + \bmh{r}_{1\tau}} \quad \text{and} \quad \theta_{12} = \dfrac{N^{0}_{1v}\bmh{\mu}_{2v}}{N^{0}_{2v}\bmh{\mu}_{1v}} \, .
\]

On the other hand,  considering the compartment  $S_{2\tau}$ in Eq. \eqref{eq:orchard2_EqS2t} for orchard $2$, we get the HLB epidemic final size of orchard 2, which is given by:
\begin{flalign}
S_{2\tau}(\infty) & = S_{2\tau}(0)\exp\left( -\dfrac{\bmh{\beta}_{21}}{N_{2\tau}}\int_0^\infty I_{2v}(t)dt \right) \cdot\exp\left(\dfrac{\bmh{\beta}_{21}p_{12}}{N_{2\tau}}\int_0^\infty I_{1v}(t)dt\right) \, ;
\end{flalign}

By following the same procedure as we did to derive Eq. \eqref{eq:Int_I1v_Final}, \emph{i.e.}, by solving $A_{2\tau} + I_{2\tau}$ in Eq. \eqref{eq:orchard2_EqI2v} and made the assumption that $I_{2v}/N_{2v} < 1$, allowing us to perform a series expansion of $1/(1-I_{2v}/N_{2v})$, we obtain the following expression:
\begin{equation}\label{eq:Int_I2v}
    \int_0^\infty I_{2v}(t)dt =  \dfrac{\bmh{\beta}_{2v} N_{2v}}{\bmh{\mu}_{2v} N_{2\tau}} \int_0^\infty (A_{2\tau} + I_{2\tau})dt  \,.
\end{equation}

And using \eqref{eq:sum_in_A1_I1} (but switching $1 \rightarrow 2$) we get:
\begin{flalign}\label{Final_Size_Trees_Orchard_2}
S_{2\tau}(\infty) & = S_{2\tau}(0)\exp\left(-\dfrac{R^{2}_{20}}{N_{2\tau}}\cdot \Biggl( R_{2\tau}(\infty) + \delta_{2} I_{2\tau}(0) \Biggr) \Biggl( 1 + \theta_{12} p^{2}_{12}\Biggr) \right) \times  \nonumber \\
& \hspace{2cm} \times \exp \left( -\dfrac{R^{2}_{10}\bmh{\beta}_{21}}{N_{2\tau}\bmh{\beta}_{11}}p_{11}p_{12}\Biggl( R_{1\tau}(\infty) + \delta_{1} I_{1\tau}(0) \Biggr) \right) \, .
\end{flalign}

The derived equations, \eqref{Final_Size_Trees_Orchard_1}-\eqref{Final_Size_Trees_Orchard_2}, offer a quantitative insight into the final size of the two-orchard model, incorporating the local basic reproduction number and the ACP emigration parameter. Furthermore, it is noteworthy that when $\phi_{12} = 0$, these equations simplify to represent the final size of each isolated orchard individually.


\begin{thebibliography}{99}

\bibitem{Li_2020}
Li, S., Wu, F., Duan, Y., Singerman, A., Guan, Z. (2020). Citrus Greening: Management Strategies and Their Economic Impact. HortScience, 55 (5), 604–612.

\bibitem{Singerman_2020}
Singerman, A., Rogers, M. E. (2020). The Economic Challenges of Dealing with Citrus Greening: The Case of Florida. Journal of Integrated Pest Management, 11 (1). 

\bibitem{Hodges_2012}
Hodges, A. W., Spreen, T. H. (2012). Economic Impacts of Citrus Greening (HLB) in Florida, 2006/07–2010/11. EDIS, 2012 (1).

\bibitem{Gottwald_2010}
Gottwald, T. R. (2010). Current Epidemiological Understanding of Citrus Huanglongbing. Annual Review of Phytopathology, 48(1), 119–139.

\bibitem{Wang_2013}
Wang, N., Trivedi, P. (2013). Citrus Huanglongbing: A Newly Relevant Disease Presents Unprecedented Challenges. Phytopathology, 103 (7), 652–665.

\bibitem{Boina_2015}
Boina, D. R., Bloomquist, J. R. (2015). Chemical control of the Asian citrus psyllid and of huanglongbing disease in citrus. Pest Management Science, 71(6), 808–823.

\bibitem{Milosavljevic_2017}
Milosavljević I., Schall K, Hoddle C., Morgan D., Hoddle M. (2017). Biocontrol program targets Asian citrus psyllid in California's urban areas. Calif. Agr. 71(3), 169-177.

\bibitem{Alquezar_2022}
Alquézar, B., Carmona, L., Bennici, S., Miranda, M. P., Bassanezi, R. B., Peña, L. (2022). Cultural Management of Huanglongbing: Current Status and Ongoing Research. Phytopathology, 112(1), 11–25.

\bibitem{Bassanezi_2020}
Bassanezi, R. B., Lopes, S. A., de Miranda, M. P., Wulff, N. A., Volpe, H. X. L., Ayres, A. J. (2020). Overview of citrus huanglongbing spread and management strategies in Brazil. Tropical Plant Pathology, 45(3), 251–264.

\bibitem{Chiyaka_2012}
Chiyaka, C., Singer, B. H., Halbert, S. E., Morris, J. G., van Bruggen, A. H. C. (2012). Modeling huanglongbing transmission within a citrus tree. Proceedings of the National Academy of Sciences, 109(30), 12213–12218.

\bibitem{Lee_2015}
Lee, J. A., Halbert, S. E., Dawson, W. O., Robertson, C. J., Keesling, J. E., Singer, B. H. (2015). Asymptomatic spread of Huanglongbing and implications for disease control. Proceedings of the National Academy of Sciences, 112(24), 7605–7610.

\bibitem{Ling_2021}
Ling, W., Wu, P., Li, X., Xie, L. (2021). Huanglongbing Model under the Control Strategy of Discontinuous Removal of Infected Trees. Symmetry, 13(7), 1164.

\bibitem{FZhang_2020}
Zhang, F., Qiu, Z., Zhong, B., Feng, T.,  Huang, A. (2020). Modeling Citrus Huanglongbing transmission within an orchard and its optimal control. Mathematical Biosciences and Engineering, 17(3), 2048–2069.

\bibitem{FZhang_2021}
Zhang, F., Qiu, Z., Huang, A., Zhao, X. (2021). Optimal control and cost-effectiveness analysis of a Huanglongbing model with comprehensive interventions. Applied Mathematical Modelling, 90, 719–741.

\bibitem{Taylor_2016}
Taylor, R. A., Mordecai, E. A., Gilligan, C. A., Rohr, J. R., Johnson, L. R. (2016). Mathematical models are a powerful method to understand and control the spread of Huanglongbing. PeerJ, 4, e2642.

\bibitem{Vilamiu_2012}
Vilamiu, R. G. d’A., Ternes, S., Braga, G. A.,  Laranjeira, F. F. (2012). A model for Huanglongbing spread between citrus plants including delay times and human intervention. AIP Conference Proceedings.

\bibitem{Lewis_Rosenblum_2015}
Lewis-Rosenblum, H., Martini, X., Tiwari, S., Stelinski, L. L. (2015). Seasonal Movement Patterns and Long-Range Dispersal of Asian Citrus Psyllid in Florida Citrus. Journal of Economic Entomology, 108(1), 3–10.

\bibitem{Tiwari_2010}
Tiwari, S., Lewis-Rosenblum, H., Pelz-Stelinski, K.,  Stelinski, L. L. (2010). Incidence of Candidatus Liberibacter asiaticus Infection in Abandoned Citrus Occurring in Proximity to Commercially Managed Groves. Journal of Economic Entomology, 103(6), 1972–1978.

\bibitem{Martini_2013}
Martini, X., Addison, T., Fleming, B., Jackson, I., Pelz-Stelinski, K., Stelinski, L. L. (2013). Occurrence of Diaphorina citri (Hemiptera: Liviidae) in an Unexpected Ecosystem: The Lake Kissimmee State Park Forest, Florida. Florida Entomologist, 96(2), 658–660.

\bibitem{Stelinski_2019}
Stelinski. (2019). Ecological Aspects of the Vector-Borne Bacterial Disease, Citrus Greening (Huanglongbing): Dispersal and Host Use by Asian Citrus Psyllid, Diaphorina Citri Kuwayama. Insects, 10(7), 208.

\bibitem{Tomaseto_2017}
Tomaseto, A. F., Miranda, M. P., Moral, R. A., de Lara, I. A. R., Fereres, A., Lopes, J. R. S. (2017). Environmental conditions for Diaphorina citri Kuwayama (Hemiptera: Liviidae) take-off. Journal of Applied Entomology, 142(1–2), 104–113.

\bibitem{Zorzenon_2020}
Zorzenon, F. P. F., Tomaseto, A. F., Daugherty, M. P., Lopes, J. R. S.,  Miranda, M. P. (2020). Factors associated with Diaphorina citri immigration into commercial citrus orchards in São Paulo State, Brazil. Journal of Applied Entomology, 145(4), 326–335. 

\bibitem{Hall_2011}
Hall, D. G., Hentz, M. G. (2011). Seasonal flight activity by the Asian citrus psyllid in east central Florida. Entomologia Experimentalis et Applicata, 139(1), 75–85. 

\bibitem{Johnston_2019}
Johnston,  N., Stelinski, L. L.,  Stansly, P. (2019).  Dispersal Patterns of Diaphorina citri (Kuwayama) (Hemiptera: Liviidae) as Influenced by Citrus Grove Management and Abiotic Factors. Florida Entomologist, 102(1), 168-173.

\bibitem{Uvencio_2022}
Giménez-Mujica, U. J., Anzo-Hernández, A., Velázquez-Castro, J. (2022). Epidemic local final size in a metapopulation network as indicator of geographical priority for control strategies in SIR type diseases. Mathematical Biosciences, 343, 108730.

\bibitem{Anzo_2019}
Anzo-Hernández, A., Bonilla-Capilla, B., Velázquez-Castro, J., Soto-Bajo, M., Fraguela-Collar, A. (2019). The risk matrix of vector-borne diseases in metapopulation networks and its relation with local and global R0. Communications in Nonlinear Science and Numerical Simulation, 68, 1–14.

\bibitem{Luo_2017}
Luo, L., Gao, S., Ge, Y., Luo, Y. (2017). Transmission dynamics of a Huanglongbing model with cross protection. Advances in Difference Equations, 2017(1).

\bibitem{Liu_2000}
Liu, Y. H.,  Tsai, J. H. (2000). Effects of temperature on biology and life table parameters of the Asian citrus psyllid, Diaphorina citri Kuwayama (Homoptera: Psyllidae). Annals of Applied Biology, 137(3), 201–206.

\bibitem{Sule_2012}
Sule, H., Muhamad, R., Omar, D., Hee, A. (2012). Life Table and Demographic Parameters of Asian Citrus Psyllid Diaphorina citri on Limau Madu Citrus suhuiensis. Journal of Entomology, 9(3), 146–154.

\bibitem{Wang_2014}
Wang, J., Gao, S., Luo, Y., Xie, D. (2014). Threshold dynamics of a Huanglongbing model with logistic growth in periodic environments. Abstract and Applied Analysis, 2014, 1–10.

\bibitem{Diekmann_2010}
Diekmann, O., Heesterbeek, J. A. P., Roberts, M. G. (2009). The construction of next-generation matrices for compartmental epidemic models. Journal of The Royal Society Interface, 7(47), 873–885.

\bibitem{Blower_1994}
Blower, S. M., Dowlatabadi, H. (1994). Sensitivity and Uncertainty Analysis of Complex Models of Disease Transmission: An HIV Model, as an Example. International Statistical Review / Revue Internationale de Statistique, 62(2), 229.

\bibitem{Seaholm_1988}
Seaholm, S. K., Ackerman, E., Wu, S.-C. (1988). Latin hypercube sampling and the sensitivity analysis of a Monte Carlo epidemic model. International Journal of Bio-Medical Computing, 23(1–2), 97–112.

\bibitem{Helton_2005}
Helton, J. C., Davis, F. J., Johnson, J. D. (2005). A comparison of uncertainty and sensitivity analysis results obtained with random and Latin hypercube sampling. Reliability Engineering amp; System Safety, 89(3), 305–330.

\bibitem{Bidah_2020}
Bidah, S., Zakary, O., Rachik, M. (2020). Stability and Global Sensitivity Analysis for an Agree-Disagree Model: Partial Rank Correlation Coefficient and Latin Hypercube Sampling Methods. International Journal of Differential Equations, 2020, 1–14.

\bibitem{Marino_2008}
Marino, S., Hogue, I. B., Ray, C. J., Kirschner, D. E. (2008). A methodology for performing global uncertainty and sensitivity analysis in systems biology. Journal of Theoretical Biology, 254(1), 178–196.

\bibitem{pyDOE}
Baudin, M. (2021). {pyDOE: The Python Design of Experiments Library}. \url{https://pypi.org/project/pyDOE/}. [Online; accessed May 05, 2023].

\bibitem{Holland_1992}
Holland, J. H. (1992). Adaptation in Natural and Artificial Systems (2nd ed.). University of Michigan Press.

\bibitem{Mitchell_1998}
Mitchell, M. (1998). An Introduction to Genetic Algorithms (1st ed.). MIT Press.
\end{thebibliography}
\end{document}